\documentclass[12pt]{article}

\usepackage[T2A]{fontenc}
\usepackage[utf8]{inputenc}

\usepackage{amsfonts}
\usepackage{amssymb}
\usepackage{amsmath}
\usepackage{amsthm}

\def \le {\leqslant}
\def \ge {\geqslant}

\theoremstyle{plain}

\usepackage{graphicx}
\usepackage{xcolor}

\begin{document}
\begin{Huge}
 \centerline{On badly approximable  vectors}
\end{Huge}
\vskip+0.5cm
\centerline{by  {\bf  Renat Akhunzhanov} 
 and {\bf Nikolay  Moshchevitin}\footnote{Research is supported by the Russian Science Foundation under grant 19-11-00001.
 }}
\vskip+1cm

Motivated by a wonderful paper \cite{Roy} where a powerful method was introduced,
we prove a criterion for a vector
$\pmb{\alpha}\in \mathbb{R}^d$ to  be a badly approximable vector. Moreover we construct certain examples which show
that a more general version of our criterion is not valid.

  \vskip+0.3cm
{\bf 1. Badly approximable real numbers and continued fractions.}
  \vskip+0.3cm

Let $||x|| = \min_{a\in \mathbb{Z}}|x-a|$ denote the distance from a real $x$  to the nearest integer.
A real  irrational  number   $\alpha$ is  called {\it badly approximable} if 
$$
\inf_{q\in \mathbb{Z}_+}  q \, ||q\alpha ||  >0.
$$
It is a well known fact that $\alpha $ is a badly approximable number if and only if 
the partial quotients 
in  continued fraction expansion
\begin{equation}\label{1m}
   [a_0;a_1,a_2,...,a_\nu,...]
   =   a_0  +
\frac{1}{\displaystyle{a_1+\frac{1}{\displaystyle{a_2 + \cdots+
\frac{1}{\displaystyle{a_\nu + ...
{} }}}}}}
   , a_0\in \mathbb{Z},\,\,\,\ a_j \in \mathbb{Z}_+, j=1,2,3,... 
\end{equation}
 are bounded, that is
$$
\sup_{\nu \ge 1} a_\nu < \infty
$$
(see, for example Theorem 5F from  Chapter I  from \cite{SCH}).
 Let us consider the sequence of  the best approximations to  $\alpha$, that is the sequence of integers
$$
q_1< q_2 <...<q_{\nu}<q_{\nu+1}<...
$$
such that 
$$
||q_\nu\alpha|| = |q_\nu\alpha - p_\nu| < ||q\alpha || ,\,\,\,\,\,\,
\text{for all positive integers } \,\,\,\,\,\,\, q< q_\nu.
$$
By Lagrange's theorem all the best approximations $(q_\nu,p_\nu)$ with $ q_\nu >1$  are just the convergents
$$\frac{p_\nu}{q_{\nu}} = [a_0; a_1,...,a_\nu]$$
  for  the  continued fraction expansion (\ref{1m}).
  For the convergents' denominators and for the remainders $ \xi_\nu = ||q_\nu\alpha||$ we have recurrent formulas
  $$
  q_{\nu+1} = a_{\nu+1} q_\nu + q_{\nu-1},\,\,\,\,\,
   \xi_{\nu+1} =  \xi_{\nu-1} -     a_{\nu+1}  \xi_\nu.
   $$
   So by taking  integer parts we have
   $$
   a_{\nu+1} = \left[\frac{q_{\nu+1}}{q_\nu}\right] = 
   \left[\frac{\xi_{\nu-1}}{\xi_\nu}\right],
   $$
   and the following obvious statement is valid.
   
   {\bf Proposition 1.}\,
   {\it An irrational number $\alpha$ is badly approximable if and only if
   $$
   \sup_{\nu \ge 1} \frac{q_{\nu+1}}{q_\nu} <\infty
   $$
   and if and only if
   $$
   \inf_{\nu \ge 1} \frac{\xi_{\nu+1}}{\xi_{\nu}} >0.
   $$
   }

In the present paper we deal with a generalization of Proposition 1 to simultaneous Diophantine approximation for several real numbers and to Diophantine approximation for one linear form.
In the next section we recall all the necessary definitions and in Section 3 
we formulate our main results.

  \vskip+0.3cm
{\bf 2. Simultaneous approximation to $d$ numbers and linear forms.}

  \vskip+0.3cm

We consider a real vector  $  \pmb{\alpha} = (\alpha_1,...,\alpha_d)\in\mathbb{R}^d$  such that  $1,\alpha_1,...,\alpha_d$ are linearly  independent over $\mathbb{Z}$.
Vector $\pmb{\alpha}$ is called {\it badly approximable} if 
\begin{equation}\label{d1}
 \inf_{q\in \mathbb{Z}_+} q^{1/d} \max_{1\le j \le d} ||q\alpha_j|| >0
 \end{equation}
 By the famous Perron-Khintchine's  transference theorem  (see Theorem 5B from Chapter IV from \cite{SCH}) condition (\ref{d1}) is equivalent to
\begin{equation}\label{d2}
 \inf_{\pmb{m} = (m_1,...,m_d) \in \mathbb{Z}^d \setminus\{\pmb{0}\}} \left(\max_{1\le j \le d} |m_j|\right)^d ||m_1\alpha_1+...+m_d\alpha_d || >0.
 \end{equation}
We  consider the best approximation vectors  for simultaneous approximation
 \begin{equation}\label{zednu}
\pmb{z}_\nu = (q_\nu, a_{1,\nu},...,a_{d,\nu}),\,\,\,\,\, \nu = 1,2,3,... \,\,\,,
\end{equation}
satisfying
$$
q_1< q_2 <...<q_{\nu}<q_{\nu+1}<...\,\,\, ,
$$
$$
\xi_\nu = \max_{1\le j \le d} ||q_\nu\alpha_j|| =   \max_{1\le j \le d} |q_\nu\alpha_j - a_{j,\nu}| 
<  \max_{1\le j \le d} ||q\alpha_j|| , \,\,\,\,\,\, \forall{q} < q_\nu,
$$
\begin{equation}\label{00}
\xi_1>\xi_2>...>\xi_{\nu}>\xi_{\nu+1}>... \,\,\,,
\end{equation}
as well as the best approximation vectors  in the sense of the linear form
 \begin{equation}\label{form}
\pmb{m}_\nu = (m_{0,\nu}, m_{1,\nu},...,m_{d,\nu}),\,\,\,\,\, \nu = 1,2,3,... \,\,\,.
\end{equation}
Namely, if we define
$M_\nu = 
\max_{1\le j\le d}  |m_{j,\nu}|$, we have
\begin{equation}\label{mioX} 
M_1<M_2<...<M_\nu<M_{\nu+1}<... \,\, .
\end{equation}
At the same time for the values of linear form
$$
L_\nu =
||m_{1,\nu}\alpha_1+...+m_{d,\nu}\alpha_d ||
=
|m_{0,\nu}+m_{1,\nu}\alpha_1+...+m_{d,\nu}\alpha_d| 
$$
the inequalities
$$
L_\nu<
||m_{1}\alpha_1+...+m_{d}\alpha_d||,\,\,\,\,\,\,
\forall (m_1,...,m_d) \in \mathbb{Z}^d \setminus\{\pmb{0}\}\,\,\,
\text{with}\,\,\,
\max_{1\le j \le d} |m_j|< M_\nu,
$$ 
and
$$
L_1>L_2>...>L_\nu>L_{\nu+1}>...   
$$ 
are valid.
Basic facts about best approximation vectors can be found for example in \cite{Che} and \cite{msing}.
In particular, from the Minkowski convex body theorem it follows that
 \begin{equation}\label{bo} 
\xi_\nu \le \frac{1}{q_{\nu+1}^{1/d}}
 \end{equation}
 and
  \begin{equation}\label{bo1}
 L_\nu \le \frac{1}{M_{\nu+1}^{d}}
 \end{equation}

    \vskip+0.3cm
  {\bf 3. Main results.}
  
    \vskip+0.3cm
    
    Our first result is the following criterium of badly approximability.

  \vskip+0.3cm
 {\bf Theorem 1.}
 \,
 {\it Suppose that $\alpha_1,...,\alpha_d, 1$ are linearly independent over $\mathbb{Q}$. Then  the following three statements are equivalent:
 
 {\rm ({\bf i})} 
 $\pmb{\alpha}$ is badly approximable;

  {\rm ({\bf ii})}  $\sup_{j} \frac{q_{j+1}}{q_j} < \infty$;

  {\rm ({\bf iii})}  $\inf_{j} \frac{L_{j+1}}{L_j} >0 $. 
  }
  \vskip+0.3cm
  
   We prove the  implication
 ({\bf ii})$\Longrightarrow$({\bf i}) 
 in Sections 6, 7.  A proof of  the implication 
 ({\bf iii})$\Longrightarrow$({\bf i}) 
 will be given in Section 8.
  Here we should note that the
 implications ({\bf i})$\Longrightarrow$({\bf ii}) and  ({\bf i})$\Longrightarrow$({\bf iii}) are obvious.
  Indeed from the definition  (\ref{d1}) and  inequality  (\ref{bo}) we immediately get
  $$
  \frac{\gamma}{q_\nu^{1/d}} \le \xi_\nu \le \frac{1}{q_{\nu+1}^{1/d}}\,\,\,\,\,\,\forall \nu
  $$
  for some positive $\gamma$ and so $ \frac{q_{\nu+1}}{q_\nu} \le \gamma^{-d},$ that is ({\bf ii}). Similarly from   (\ref{d2})   we get
  $$
   L_{\nu+1} \ge \frac{\gamma}{M_{\nu+1}^d} \,\,\,\,\,\, \forall \nu
  $$
  with some positive $\gamma$ and together with 
   (\ref{bo1})  this gives
   $$
   \frac{L_{\nu+1}}{L_\nu} \ge \gamma,
   $$ 
  and this is ({\bf iii}).

 In fact  for badly approximable $\pmb{\alpha}$ we can say something more, by the same argument.
   
    \vskip+0.3cm
  {\bf Remark 1.}
 \,
 {\it  If $\pmb{\alpha}\in \mathbb{R}^d$ is badly approximable then besides the inequalities {\rm  ({\bf ii})} and  {\rm ({\bf iii})} 
 the inequalities 
 \begin{equation}\label{dg}
 \inf_{j} \frac{\xi_{j+1}}{\xi_j} >0
 ,\,\,\,\,\,\text{and}\,\,\,\,\,\,
  \sup_{j} \frac{M_{j+1}}{M_j} < \infty
 \end{equation}
 are also valid.} 
    
    \vskip+0.3cm
    Indeed, we can easily get the first inequality from (\ref{dg}) by combining inequality
    $ \xi_{\nu+1} >\gamma/q_{\nu+1}^{1/d}$ and (\ref{bo}); the second inequality from (\ref{dg}) can be obtained by combining
    $L_\nu > \gamma M_\nu^{-d}$ and (\ref{bo1}). However the converse statements are not true.
    Our second result is  given by the following statement. For  the simplicity reason we formulate and prove this result for two-dimensional case only.
    However the construction may be easily generalized to  the case of simultaneous approximation to $d$ numbers.

   \vskip+0.3cm
 {\bf Theorem 2.}
 \,
 {\it  There exists uncountably many $\pmb{\alpha} = (\alpha_1,\alpha_2)\in \mathbb{R}^2$ such that 
 
 $\bullet$\,
 $1,\alpha_1,\alpha_2$ are linearly independent over $\mathbb{Z}$;
 
 $\bullet$\,
 $\inf_\nu\frac {\xi_{\nu+1}}{\xi_\nu} >0$;
 
 $\bullet$\, $\pmb{\alpha} $ is not badly approximable.
 }
 
   \vskip+0.3cm
   
   The idea of the construction from the proof of Theorem 2 is quite simple. 
   It is related to a construction from our earlier paper \cite{mou}.
   One should construct a vector $\pmb{\alpha} \in \mathbb{R}^2$
    such that the best approximation vectors to it  for long times lie in two-dimensional subspaces. Moreover, for the integer approximations from these two-dimensional subspaces we should  ensure  some kind of "one-dimensional badly approximability". However a complete proof  for Theorem 2 is rather cumbersome.    We give our proof of Theorem 2 in Sections 9, 10 and 11.
    
    We would like to note that very recently during the refereeing process of this paper an alternative construction to  prove Theorem 2 by means of Parametric Geometry of Numbers  based on  on a deep theorem due to D. Roy 
    \cite{R} was obtained by W.M. Schmidt \cite{SM}.
    
     \vskip+0.3cm
     
   In the present paper we would like to announce a theorem dual to Theorem 2 which deals with the best approximations in the sense of a linear form.  
   The formulation of this result is below.
   
      \vskip+0.3cm
 {\bf Theorem 3.}
 \,
 {\it  There exist uncountably many $\pmb{\alpha} = (\alpha_1,\alpha_2)\in \mathbb{R}^2$ such that 
 
 $\bullet$\,
 $1,\alpha_1,\alpha_2$ are linearly independent over $\mathbb{Z}$;
 
 $\bullet$\,
 $\sup_\nu\frac {M_{\nu+1}}{M_\nu} < \infty$;
 
 $\bullet$\, $\pmb{\alpha} $ is not badly approximable.
 }
 
   \vskip+0.3cm
   
   In this paper we do not give a proof of Theorem 3 but just announce it. The proof we have is based on the same idea as the proof of Theorem 2 but  it is even more  technical and  cumbersome. 
   Moreover, it is related  to some general phenomenon, and we suppose to consider it in a separate paper which now is in preparation.
   


       \vskip+0.3cm
 
 {\bf  4. On Diophantine exponents.}
 
   \vskip+0.3cm
   
 For a real  $\pmb{\alpha} \in \mathbb{R}^d$ we recall the definitions of the {\it 
 ordinary } Diophantine exponent $\omega(\pmb{\alpha}) $ and the  {\it uniform} Diophantine exponent
  $\hat{\omega}(\pmb{\alpha}) $  in the sense of  simultaneous Diophantine approximation.
  The ordinary  Diophantine exponent  $\omega(\pmb{\alpha}) $ is defined  as the supremum of those 
  $\gamma \in \mathbb{R}$ for which there exists  an unbounded sequence of values of  $T\in \mathbb{R}_+$ such that the system 
  \begin{equation}\label{Xdiop}
  \begin{cases}
  \displaystyle{
  \max_{1\le j \le d}}
  ||q\alpha_j || \le T^{-\gamma},
  \cr
  1\le q \le T
  \end{cases}
  \end{equation}
  has an integer solution $ q\in \mathbb{Z}$.   The uniform Diophantine exponent  $\hat{\omega}(\pmb{\alpha}) $ is defined  as the supremum of those 
  $\gamma \in \mathbb{R}$ for which there exists  $T_0$ such that for every $ T\ge T_0$  the system  (\ref{Xdiop})
  has an integer solution $ q\in \mathbb{Z}$. 
   Equivalently in terms of the best approximation vectors, $\hat{\omega}(\pmb{\alpha}) $ can be defined as the supremum of those 
    $\gamma \in \mathbb{R}$ for which  the inequality
    \begin{equation}\label{uuw}
    \xi_\nu\le q_{\nu+1}^{-\gamma}
    \end{equation}
    is  valid for all $\nu$ large enough.

  It is well known that
  $$
  \frac{1}{d}\le \hat{\omega}(\pmb{\alpha})  \le 1
  $$
  for every $\pmb{\alpha} \in \mathbb{R}^d\setminus \mathbb{Q}^d$ and obviously
  $$
  \hat{\omega}(\pmb{\alpha}) \le {\omega}(\pmb{\alpha}) \le +\infty .
  $$

  As it was discovered by  V. Jarn\'{\i}k \cite{JRUS},
  the first  trivial inequality  here can be improved. The optimal lowed bound for 
  ${\omega}(\pmb{\alpha}) $ in terms of  $\hat{\omega}(\pmb{\alpha}) $ was obtained in  \cite{MaMo}
  where the authors solve a problem by W.M. Schmidt and L. Summerer \cite{SS}.   
            In the case when the numbers $1, \alpha_1,...,\alpha_d$ are linearly independent over $\mathbb{Q}$ in the paper  \cite{MaMo}
            the authors establish the inequality
       \begin{equation}\label{moma}
       \frac{\omega (\pmb{\alpha})}{
            \hat{\omega} (\pmb{\alpha}) }
            \ge G_d (\hat{\omega} (\pmb{\alpha})),
            \end{equation}
            where $ G_d (\hat{\omega} (\pmb{\alpha}))\ge 1$ is the positive root of the equation
             \begin{equation}\label{momaeq}
      t^{d-1} =  \frac{  \hat{\omega} (\pmb{\alpha})}{1-  \hat{\omega} (\pmb{\alpha})} (1+t+...+t^{d-2}).
            \end{equation}
            The main argument of the proof from  \cite{MaMo} is that there exist infinitely many $\nu$  with
   \begin{equation}\label{momaqu}
   q_{\nu+1} \ge q_\nu^{ G_d (\hat{\omega} (\pmb{\alpha}))}.
   \end{equation}

 Here we should note that the wonderful paper \cite{Roy}  
  deals with a simple and elegant proof of this result as well as with some other  related problems.
  
  If $\pmb{\alpha} \in \mathbb{R}^d$ is a badly approximable vector we have  $ {\omega}(\pmb{\alpha}) = \hat{\omega}(\pmb{\alpha}) = \frac{1}{d}$.
  
  However, Theorem 2 shows that for $ d \ge 2$ the condition 
  \begin{equation}\label{uzas}
  \inf_\nu\frac {\xi_{\nu+1}}{\xi_\nu} >0
  \end{equation}
 may be satisfied for $\pmb{\alpha}$ which is not badly approximable. Moreover the construction from the proof of Theorem 2 gives $\pmb{\alpha}$ with
 $ \hat{\omega}(\pmb{\alpha}) = \frac{1}{2}$ and $ {\omega}(\pmb{\alpha}) =1$. We would like to give a comment on this, and formulate the following statement.
 
    \vskip+0.3cm
 {\bf Proposition 2.}\, {\it
 Suppose  that among the numbers $ \alpha_1, ...,\alpha_d$ there  exist at least two numbers linearly independent together with $1$ over $\mathbb{Q}$, and suppose that  $\pmb{\alpha}$ satisfies condition (\ref{uzas}).
 Then 
 \begin{equation}\label{een}
   \hat{\omega}(\pmb{\alpha}) \le \frac{1}{2}.
   \end{equation}
 }
   
      \vskip+0.3cm
      Proof.  \,  Jarn\'{\i}k \cite{JRUS} proved that under the conditions of Proposition 1 there exist infinitely many  linearly independent triples $\pmb{z}_{\nu-1},\pmb{z}_\nu, \pmb{z}_{\nu+1}$  of consecutive best approximation vectors. Moreover   
      for such a triple there exist indices $ j_1,j_2$ such that 
      $$
      D =
                \left|
      \begin{array}{ccc}
      q_{\nu-1} & a_{j_1,\nu-1}& a_{j_2, \nu-1}
      \cr
      q_{\nu} & a_{j_1,\nu}& a_{j_2, \nu}
      \cr
      q_{\nu+1} & a_{j_1,\nu+1}& a_{j_2, \nu+1}
      \end{array}
      \right|
      =
      \left|
      \begin{array}{ccc}
      q_{\nu-1} & a_{j_1,\nu-1}-  q_{\nu-1}\alpha_{j_1}& a_{j_2, \nu-1} -  q_{\nu-1}\alpha_{j_2}
      \cr
      q_{\nu} & a_{j_1,\nu} -  q_{\nu}\alpha_{j_1}& a_{j_2, \nu} -  q_{\nu}\alpha_{j_2}
      \cr
      q_{\nu+1} & a_{j_1,\nu+1} -  q_{\nu+1}\alpha_{j_1}& a_{j_2, \nu+1} -  q_{\nu+1}\alpha_{j_2}
      \end{array}
      \right| \neq 0.
      $$
      But from the definition of values $\xi_\nu$ and (\ref{uzas}) we see that 
      $$
      1\le |D| \le
      6 \xi_{\nu-1}\xi_\nu q_{\nu+1}
      \ll
      \xi_\nu^2 q_{\nu+1}
      $$
      (of course here the constant in the sign $\ll$ may depend on $\pmb{\alpha}$).  The last inequality together with the definition of $ \hat{\omega}(\pmb{\alpha}) $ in terms of the inequality   (\ref{uuw})
      gives (\ref{een}).$\Box$
      
         \vskip+0.3cm
         It is clear that the bound (\ref{een}) is optimal for $d=2$. However what are  admissible values of $\hat{\omega}(\pmb{\alpha})$ and  ${\omega}(\pmb{\alpha})$ 
         for general $d$  
         under the condition (\ref{uzas})
         for the  numbers $1, \alpha_1,...,\alpha_d$  which are linearly independent over $\mathbb{Q}$ seems to be an open question.

               \vskip+0.3cm
            In addition, here we would like to give the following 
            remark.  
            We should note that if
 \begin{equation}\label{roy1}
 bq_\nu^{-\beta} < \xi_\nu < aq_{\nu+1}^{-\alpha}
 \end{equation}
 with some positive $a,b$ and $\alpha,\beta$  satisfying $ \beta\ge\alpha \ge 1/d$, then
\begin{equation}\label{wu?}
 q_{\nu+1} < Cq_\nu^{  {\beta}/{\alpha} }\,\,\,\,\,\text{with}\,\,\,\,C = \left(\frac{a}{b}\right)^{{1}/{\alpha}},
 \end{equation}
 in particular
 $$
  q_{\nu+1} < C'q_\nu^{d{\beta}}\,\,\,\,\,\text{with}\,\,\,\,\, C' = \frac{1}{b^d}.
 $$

            Consider the exponent
            $$
            \tau(\pmb{\alpha}) = \limsup_{\nu\to \infty} \frac{\log q_{\nu+1}}{\log q_\nu} 
            $$
            which contain information about the growth of the best approximation vectors to $\pmb{\alpha}$. Then the observation mentioned above can be summarized as
             \vskip+0.3cm
 {\bf Proposition 3.}\, {\it  Suppose that  the numbers $1, \alpha_1,...,\alpha_d$ are linearly independent over $\mathbb{Q}$. Then
 \begin{equation}\label{priu1}
 G_d (\hat{\omega} (\pmb{\alpha}) )\le   \tau(\pmb{\alpha}) \le \frac{\omega(\pmb{\alpha})}{\hat{\omega}(\pmb{\alpha})} \le d\omega(\pmb{\alpha}).
 \end{equation}
 Moreover
  \begin{equation}\label{priu2}
 \hat{\omega}(\pmb{\alpha}) \le \frac{1}{\sum_{j=0}^{d-1}   \tau(\pmb{\alpha})^{-j}}.
 \end{equation}

          }
                     \vskip+0.3cm    
                     
                     Proof. Lower bound  for 
                     $\tau(\pmb{\alpha})$ in 
                      (\ref{priu1}) immediately follows from (\ref{momaqu}). Upper bound comes from (\ref{wu?}) under the condition (\ref{roy1}).
                     Inequality  (\ref{priu2})  follows from (\ref{momaqu}) and (\ref{momaeq}).$\Box$
               
             \vskip+0.3cm
 
 {\bf  5. Some notation.}
 
    \vskip+0.3cm

 We use the following notation.  
 Together with the best approximation vectors (\ref{zednu}) which we have denoted by $\pmb{z}_\nu$ 
 we consider the points
\begin{equation}\label{ze}
\pmb{Z}_\nu = (q_\nu, q_\nu\alpha_{1},..., q_\nu \alpha_{d}).
\end{equation}
By 
$|\pmb{\xi}|
$ we denote the Euclidean norm of the vector $\pmb{\xi}  \in \mathbb{R}^{k}$ in any dimension $k$.
By 
$$
|\pmb{\eta}|_\infty
=\max_{1\le j \le d} |\eta_j|
$$
we denote the sup-norm of the vector $\pmb{\eta} \in \mathbb{R}^d$.
 In the case $ \pmb{x}  = (x_0,x_1,...,x_d) \in \mathbb{R}^{d+1}$
we will use the notation
$$
    |\underline{\pmb{x}}|_\infty =
\max_{1\le j \le d} |x_j|
$$
to deal with the sup-norm of the shortened vector $
\underline{\pmb{x}} = (x_1,...,x_d) \in \mathbb{R}^d$.
So
for
$
\pmb{\xi}_\nu =  \pmb{Z}_\nu - \pmb{z}_\nu
$
we have
$\xi_\nu  
=
|\underline{\pmb{\xi}}_\nu |_\infty.
$

 It is clear that 
\begin{equation}\label{0}
\pmb{\xi}_\nu   = | \pmb{Z}_\nu - \pmb{z}_\nu|\le \sqrt{d} \,  \xi_\nu
 \end{equation}
 Let 
 $$\rho(\mathcal{A},\mathcal{B})
 = \inf_{\pmb{a} \in A, \, \pmb{b} \in B} |a-b|
 $$ be the Euclidean distance between sets $ \mathcal{A}$ and $\mathcal{B}$.

  \vskip+0.3cm
 {\bf  6. Main geometric lemma.}
  
    \vskip+0.3cm

We define inductively a special collection of $d+1$ 
linearly independent best approximation vectors.
Let
$
\nu_1 = \nu,\,\, \nu_2 = \nu+1.
$ Then,
if $\pmb{z}_{\nu_1}, \pmb{z}_{\nu_2},...,\pmb{z}_{\nu_{j-1}}$ are defined we find the smallest $ \mu \ge \nu_{j-1}+1$ such that
 the vectors 
$\pmb{z}_{\nu_1}, \pmb{z}_{\nu_2},...,\pmb{z}_{\nu_{j-1}}, \pmb{z}_\mu$  are independent and put
$ \pmb{z}_{\nu_{j}}=\pmb{z}_\mu$.
At the end of the procedure we have $d+1$ independent vectors
\begin{equation}\label{points}
\pmb{z}_{\nu_1}, \pmb{z}_{\nu_2},...,\pmb{z}_{\nu_{d+1}}.
\end{equation}
We define  linear subspaces 
\begin{equation}\label{planes}
\pi_j = \langle \pmb{z}_{\nu_1},\pmb{z}_{\nu_2},...,\pmb{z}_{\nu_j}\rangle_{\mathbb{R}},\,\,\,\,
j = 1,...,d+1
\end{equation}
and  lattices 
\begin{equation}\label{lattices}
\Gamma_j = \pi_j \cap \mathbb{Z}^{d+1}.
\end{equation}
In particular  $ \Gamma_1 = \langle \pmb{z}_\nu\rangle_{\mathbb{Z}}$ and $ \Gamma_{d+1} = \mathbb{Z}^{d+1}$.
By $ \Delta_j$ we denote the $j$-dimensional fundamental volume of lattice  $
\Gamma_j$. In particular $\Delta_1 = |\pmb{z}_\nu|$ and $ \Delta_{d+1} = 1$.

Here we should note that
by Minkowski Convex Body Theorem applied for the two-dimensional lattice $\Gamma_2$ we have

 \begin{equation}\label{0000}
\xi_\nu q_{\nu+1} \le {\Delta_2},
\end{equation}
and also 
 \begin{equation}\label{000}
\xi_\nu q_{\nu+1} \ge K{\Delta_2}, \,\,\,\,\,
\text{where}
\,\,\,\,\,
K= \frac{1}{2\sqrt{d(1+ \alpha_1^2+...+\alpha_d^2)}}
 \end{equation}
 (for the details see for example \cite{msing} or Theorem 1.5  from \cite{Cheung}). Moreover, (\ref{000}) together with 
 (\ref{bo}) for every best approximation $ \nu \ge 1$ 
 gives
 $$
 K\Delta_2 \le q_{\nu+1}^{\frac{d-1}{d}}
 $$
 or
 \begin{equation}\label{dee}
(K \Delta_2)^{\frac{d}{d-1}}\le  q_{\nu+1}.
 \end{equation}

\vskip+0.3cm

 {\bf Lemma 1.}\, {\it  For every $j$ one has 
 $$
 \frac{\Delta_{j+1}}{\Delta_j} 
 \le  2\sqrt{d}\,\frac{q_{\nu_{j+1}}}{q_{\nu_{j+1}-1} } \,\xi_{\nu_{j+1}-1}.
 $$
 }
 
 \vskip+0.3cm
 
 Proof.  
 Let $\pmb{w}\in \Gamma_{j+1}\setminus\Gamma_j$ be a primitive vector such that
 $$
 \Gamma_{j+1} = \langle  \Gamma_j ,\pmb{w}\rangle_{\mathbb{Z}}
 $$
 It is clear that the lattice $\Gamma_{j+1}$ splits into a union of affine sublattices with respect to $\Gamma_j$:
 $$
 \Gamma_{j+1} = \bigcup_{k\in \mathbb{Z}} \left( \Gamma_{j}+  k\pmb{w}\right).
 $$
 We consider affine $j$-dimensional subspaces 
 $$
 \pi_{j,k} = \pi_j + k\pmb{w} \supset  \Gamma_{j}+  k\pmb{w}.
 $$
 It is clear that the Euclidean distance between  each two neighboring subspaces 
 $\pi_{j,k}$ and $\pi_{j,k+1}$ is equal to $ \frac{\Delta{j+1}}{\Delta_j} $. So in the case $ k\neq 0$ we have
 \begin{equation}\label{q}
 \rho (\pi_\nu, \pi_{\nu,k} )=   |k| \cdot \frac{\Delta_{j+1}}{\Delta_j}  \ge \frac{\Delta_{j+1}}{\Delta_j}.
 \end{equation}
Define $ k_*$ from the condition
 $$
 \pmb{z}_{\nu_{j+1}} \in   \pi_{j,k_*}.
 $$
 As $  \pmb{z}_{\nu_{j+1}}  \not\in \pi_j$ we have $ k_* \neq 0$.
 As $ \pmb{z}_{\nu_{j+1}-1}\in \pi_j$ 
  from (\ref{0}) we get
  $$
 \rho (\pmb{Z}_{\nu_{j+1}-1},  \pi_j)\le \sqrt{d} \,\xi_{\nu_{j+1}-1}.
$$
As
$$
\frac{
|
\pmb{Z}_{\nu_{j+1}}|
}{
|
\pmb{Z}_{\nu_{j+1}-1}|
}
=
\frac{q_{\nu_{j+1}}}{q_{\nu_{j+1}-1}}
$$
we deduce 
\begin{equation}\label{a}
 \rho (\pmb{Z}_{\nu_{j+1}-1},  \pi_j)
 =
\frac{q_{\nu_{j+1}}}{q_{\nu_{j+1}-1}}
 \cdot
  \rho (\pmb{Z}_{\nu_{j+1}},  \pi_j)
 \le 
\frac{q_{\nu_{j+1}}}{q_{\nu_{j+1}-1}}
 \cdot
 \sqrt{d} \,\xi_{\nu_{j+1}}.
\end{equation}
As  $ \pmb{z}_{\nu_{j+1}}\in \pi_{j,k_*}$ 
we see that 
\begin{equation}\label{b}
 \rho (\pmb{Z}_{\nu_{j+1}},  \pi_{j,k_*})\le \sqrt{d} \,\xi_{\nu_{j+1}}.
\end{equation}
 From (\ref{q}), triangle inequality, formulas    (\ref{a},\ref{b})  and the inequalities $\xi_{\nu_{j+1}} < \xi_{\nu_{j+1}-1} $  and
 s $q_{\nu_{j+1}} > q_{\nu_{j+1}-1} $
 we get 
 $$
 \frac{\Delta_{j+1}}{\Delta_j}\le 
 \rho (\pi_j, \pi_{j, k_*}) \le   \rho (\pmb{Z}_{\nu_{j+1}},  \pi_{j})+ 
   \rho (\pmb{Z}_{\nu_{j+1}},  \pi_{j,k_*}) \le \sqrt{d}\, \frac{q_{\nu_{j+1}}}{q_{\nu_{j+1}-1} } \,\xi_{\nu_{j+1}-1} +   \sqrt{d} \,\xi_{\nu_{j+1}}
   \le 2\sqrt{d}\, \frac{q_{\nu_{j+1}}}{q_{\nu_{j+1}-1} } \,\xi_{\nu_{j+1}-1}. 
   $$
 Everything is proved.$\Box$

    \vskip+0.3cm
  {\bf  7. Proof of Theorem 1: simultaneous approximation.}
     \vskip+0.3cm

 Let $ \alpha_1,...,\alpha_d$ be given. 
We suppose that ({\bf ii})  is valid and deduce ({\bf i}). For a given $\nu$ from (\ref{000})   and $\Delta_{d+1} = 1$  we get the inequality
  \begin{equation}\label{r}
  \xi_\nu q_{\nu+1} \ge K \Delta_2  =
  K\cdot
  \frac{\Delta_2}{\Delta_3}
  \cdot
    \frac{\Delta_3}{\Delta_4} 
    \cdots
     \frac{\Delta_d}{\Delta_{d+1}}.
     \end{equation}
     Now  we deduce from 
     ({\bf ii}) the condition ({\bf i}).
     Lemma 1 gives
     \begin{equation}\label{l}
      \xi_\nu q_{\nu+1} 
      \ge 
       \frac{K}{ (2\sqrt{d})^{d-1}}\cdot
       \prod_{j=3}^{d+1} \frac{q_{\nu_j-1}}{q_{\nu_j}} \cdot
      \frac{1}{  \prod_{j=3}^{d+1} \xi_{\nu_j-1}}.
     \end{equation}
     As we supposed that  ({\bf ii}) is valid, there exists $M$ such that 
     $$
     \frac{q_{\nu+1}}{q_\nu} \le M \,\,\,\,\,\,\, \forall \nu.
     $$
     Moreover  from (\ref{bo}) we have $ \xi_{\nu_j-1} \le  \xi_{\nu}   \, \forall j = 3, ...,d+1$. Now we continue  with  (\ref{l})  and get
     $$
        \xi_\nu q_{\nu+1} 
      \ge 
      \frac{K}{(2\sqrt{d}M\xi_\nu)^{d-1}} .
      $$
      As $ q_{\nu+1} \le Mq_\nu$ we get
      $$
      q_\nu^{1/d} \xi_\nu \ge 
       \frac{K^{1/d}}{(2\sqrt{d})^{(d-1)/d}M}
      \,\,\,\,\,\, \forall \nu
      $$
      and  ({\bf i})  is proved.

    \vskip+0.3cm
  {\bf  8. Proof of Theorem 1: linear form.}
     \vskip+0.3cm

We suppose that ({\bf iii})  is valid and deduce ({\bf i}). We follow the same argument as in Sections  5,6, but we need to make some changes.
We use a standard  trick which reduces the problem for  linear forms to the problem for simultaneous approximation. This trick was used in \cite{MaMo}, Section 5.2.

The proof is quite similar so we will give just a sketch of a proof. 
First of all we need a generalization of Lemma 1. 
Suppose that $\Lambda$ be a full-dimensional lattice in $\mathbb{R}^{d+1}$ with coordinates $(x_0,x_1,...,x_d)$. Suppose that 



\vskip+0.3cm
 
\noindent
({\bf a})  the intersection $ \Lambda \cap \{ \pmb{x} \in \mathbb{R}^{d+1}:\,\, x_0 = 0\}$ consists just of one lattice point 
$\pmb{0} \in \Lambda$, so
every affine subspace  of the form
$ \mathcal{A}_r = \{ \pmb{x} \in \mathbb{R}^{d+1}:\,\, x_0 = r\}$
consists of at most one point from $\Lambda$, that is the cardinality fo the intersection
$\mathcal{A}_r  \cap \Lambda$ is not greater than 1 for any $ r \in \mathbb{R}$.

\vskip+0.3cm 

We consider the best  simultaneous approximations of the line 
$$\ell  =  \{ \pmb{x} \in \mathbb{R}^{d+1}:\,\, x_1=x_2=...=x_d=0\}
$$ 
by the points of the lattice $\Lambda$. Here by the {\it best approximation point} we mean a point
$ \pmb{z} = (z_0, z_1,...,z_d)\in \Lambda$ such that in the parallelepiped
$$
\Pi_{\pmb{z}} = \left\{ \pmb{z}' = (z_0',z_1',...,z_d')\in \mathbb{R}^{d+1}:\,\,\,
|z_0'|\le |z_0|,\,\,\,
|\underline{\pmb{z}}'|_\infty \le |\underline{\pmb{z}}|_\infty \right\} 
$$
there is no lattice points different from the points $\pmb{0}, \pm \pmb{z}$, that is 
$$
\Pi_{\pmb{z}} \cap \Lambda = \{ \pmb{0},\pmb{z}, -\pmb{z}\}.
$$
As the condition ({\bf a}) is satisfied,
 for any best approximation vector
$\pmb{z}_1 = (z_{0,1},z_{1,1},..., z_{d,1} )\in \Lambda$ we can consider the unique  finite or  infinite sequence  of best approximation vectors
$\pmb{z}_\nu   = (z_{0,\nu},z_{1,\nu},..., z_{d,\nu} ) \in \Lambda, \, \nu =1,2,3,...$  such that 

\vskip+0.3cm
\noindent
$\bullet$ \,\,  $0< z_{0,1} < z_{0,2}<...<  z_{0,\nu}<  z_{0,\nu+1}<...$

\vskip+0.3cm
\noindent
$\bullet$ \,\,     $ |\underline{\pmb{z}}_1|_\infty> |\underline{\pmb{z}}_2|_\infty> ...> |\underline{\pmb{z}}_\nu|_\infty>
|\underline{\pmb{z}}_{\nu+1}|_\infty>.... $.

\vskip+0.3cm

\noindent
$\bullet$ \,\,  There is no lattice points in the parallelepiped
$$
\Pi_\nu = \left\{\pmb{z}' = (z_0',z_1',...,z_d') \in \mathbb{R}^{d+1}:\,\,
|z_0'| \le z_{0,\nu+1},\,\, |\underline{\pmb{z}}'|_\infty \le  |\underline{\pmb{z}}_\nu |_\infty \right\}
$$
besides the points $\pmb{0}, \pm \pmb{z}_\nu, \pm{\pmb{z}}_{\nu+1}$:
$$
\Pi_\nu \cap \Lambda = \{
\pmb{0},  \pmb{z}_\nu,  -\pmb{z}_\nu, {\pmb{z}}_{\nu+1}, - {\pmb{z}}_{\nu+1}\}.
$$

 \vskip+0.3cm
 
 The sequence of the best approximation vectors $\pmb{z}_\nu$ is infinite if 
 there is no non-zero lattice points on the  axis $ \ell $.
 If there is a non-zero point $\pmb{z} \in \Lambda \cap \ell$ then the sequence of the best approximation vectors is finite.
 In our proof we need to consider the case when this sequence is finite. 
 We suppose that our lattice $\Lambda $ and the best approximation vector $ \pmb{z}_1$ satisfy one more condition
 
 \vskip+0.3cm
 
\noindent
({\bf b})  the sequence of the best approximation vectors $\pmb{z}_\nu, \nu \ge 1$ does not lie in a proper linear subspace of $\mathbb{R}^{d+1}$.

  \vskip+0.3cm
  
  Now  for the lattice  $\Lambda$ satisfying conditions ({\bf a}) and ({\bf b}) we are able to  define points $\pmb{z}_{\nu_j}, \nu_1 = 1< \nu_2<... <\nu_{d+1} $ from  (\ref{points})  and subspaces  $\pi_j$ from (\ref{planes}).
  In the definition of lattice $ \Gamma_j$ there will be a slight difference. Instead of (\ref{lattices}) we put
  $$
  \Gamma_j = \pi_j \cap \Lambda.
  $$
 Again by $\Delta_j$ we define the fundamental volumes of $j$-dimensional  lattices $\Gamma_j$. In particular
 \begin{equation}\label{change}
  \Delta_{d+1} = {\rm det}\, \Lambda
  .
  \end{equation}
  The inequality  (\ref{000}) transforms now into the following statement.

   \vskip+0.3cm
   
   {\bf Lemma 2.} \,\, {\it Suppose that for a certain $ \nu$ we have  
   \begin{equation}\label{0002}
    |\underline{\pmb{z}}_\nu|_\infty \cdot
     {z}_{0,\nu+1} \ge 1
     \end{equation}
     and
        \begin{equation}\label{0003}
     |\underline{\pmb{z}}_\nu|_\infty \le 1.
     \end{equation}  
     Then
   \begin{equation}\label{0001}
 |\underline{\pmb{z}}_\nu|_\infty  \cdot z_{0,\nu+1} \ge \frac{\Delta_2}{2\sqrt{2} \, d}.
 \end{equation}
 }
  
  Proof. In fact, this lemma follows from inequality (59) of Lemma 10 from \cite{MaMo}. For the sake of completeness we give here a proof.
  Consider the $2\times (d+1)$ matrix
  $$
  \mathcal{M}=
  \left(
  \begin{array}{ccccc}
  z_{0,\nu}& z_{1,\nu}& z_{2,\nu}& ...& z_{d,\nu}
  \cr
   z_{0,\nu+1}& z_{1,\nu+1}& z_{2,\nu+1}& ...& z_{d,\nu+1}
  \end{array}
  \right).
  $$
  Then $\Delta_2^2$ is just the sum of squares of all $2\times 2$ minors  
  $$
  M_{i,j} =
   \left|
  \begin{array}{cc}
  z_{i,\nu}& z_{j,\nu}
  \cr
   z_{i,\nu+1}& z_{j,\nu+1}
  \end{array}
  \right|
  $$
  of matrix $\mathcal{M}$, that is 
  $$
  \Delta_2^2 = \sum_{0\le i<j\le d}    M_{i,j} ^2.
  $$
  As $ z_{0,\nu} < z_{0,\nu+1}$ and $   |\underline{\pmb{z}}_\nu|_\infty>  |\underline{\pmb{z}}_{\nu+1}|_\infty$ we have
  $$
  |  M_{0,j}|\le 2  |\underline{\pmb{z}}_\nu|_\infty  \cdot z_{0,\nu+1} ,\,\,\,\,\, \forall\, j = 1,2,...,d.
  $$
  From (\ref{0003}) we see that
    $$
  |  M_{i,j}|\le 2   ,\,\,\,\,\, \forall\,i,  j = 1,2,...,d.
  $$
  So by (\ref{0002}) we get 
  $$
   \Delta_2^2 \le
   4d (|\underline{\pmb{z}}_\nu|_\infty  \cdot z_{0,\nu+1} )^2 + 4d^2\le
   8d^2 (|\underline{\pmb{z}}_\nu|_\infty  \cdot z_{0,\nu+1} )^2, 
  $$
  and  Lemma 2 follows.$\Box$

 \vskip+0.3cm
  
 Instead of Lemma 1 now we have the following statement.
 
 \vskip+0.3cm

 {\bf Lemma 1$'$.}\, {\it   Suppose that the lattice $\Lambda$ and the best approximation  vector $\pmb{z}_1 $ satisfy properties {\rm ({\bf a})} and {\rm  ({\bf b})} and  consider the  best approximation vectors   (\ref{points}).
 Then
 for every $j$ one has 
 $$
 \frac{\Delta_{j+1}}{\Delta_j} 
 \le  2\sqrt{d}\,\frac{z_{0,\nu_{j+1}}}{z_{0,\nu_{j+1}-1} } \,|\underline{\pmb{z}}_{\nu_{j+1}-1}|_\infty.
 $$
 }
 
 \vskip+0.3cm
 
 The proof of Lemma 1$'$ just follows the steps of the proof of Lemma 1.
 The only difference is that  instead of the points $\pmb{Z}_\nu$ defined in (\ref{ze}) which lie on the line $\langle (1,\alpha_1,...,\alpha_d)\rangle_\mathbb{R}$ one should consider the points
 $$
 \pmb{Z}^\circ_\nu = (z_{0,\nu}, 0,...,0) \in \ell.
 $$
 We left the proof to  the reader.$\Box$
 
  \vskip+0.3cm
  Now we are ready to deduce badly approximability of $\pmb{\alpha}$ from the condition ({\bf iii}).
  Let us consider best approximation vectors (\ref{form}).    
  It may happen that there exists $\nu_0$ and a proper linear subspace $ \mathcal{L}\subset \mathbb{R}^{d+1}$ 
  of dimension $3\le l={\rm dim}\, \mathcal{L} < d+1$ such that $ \pmb{m}_\nu \in \mathcal{L}$ for all  $\nu\ge \nu_0$ (see \cite{MosheDAN} for the first result in this direction and \cite{msing} and the literature therein for a survey and related results).
  But  we will show later that under condition ({\bf iii}) this is not possible.
  
   So first of all we consider the case when
   for any $\nu_0$  the best approximation vectors $\pmb{m}_\nu, \nu \ge \nu_0$ do not lay in a proper linear subspace of
  $\mathbb{R}^{d+1}$.
Suppose that  vectors
 $$
 \pmb{m}_1 , \pmb{m}_2,...., \pmb{m}_\mu
 $$
 do not 
 lay in a proper linear subspace of
  $\mathbb{R}^{d+1}$.
 
 We consider the lattice
 $$
 \Lambda_{\pmb{\alpha}} =
 \left(
 \begin{array}{ccccc}
  1&  \alpha_1&\alpha_2&...&\alpha_d
    \cr
 0&1&...&0&0
 \cr
  0&0&...&0&0
  \cr
  0&0&...&1&0
  \cr
   0&0&...&0&1
 \end{array}
 \right)
 \,\,
 \mathbb{Z}^{d+1},
 $$
 a  parameter 
$
 T >0$
 and the lattice
 $$
  \Lambda_{\pmb{\alpha}}^{[\mu]} =  \mathcal{G}
  \Lambda_{\pmb{\alpha}} ,\,\,\,\,\,\,
   \mathcal{G} =
   \left(
 \begin{array}{ccccc}
 T^{-d}&0&...&0&0
 \cr
  0&T&...&0&0
  \cr
   0&0&...&T&0
   \cr
 0&0&...&0&T
 \end{array}
 \right),\,\,\,\,\,\,\,\,
 {\rm det} \, \Lambda_{\pmb{\alpha}}^{[\mu]}  = 1.
  $$
  As $ \alpha_1,...,\alpha_d$ are linearly independent over $\mathbb{Z}$,
  the lattice 
  $
  \Lambda_{\pmb{\alpha}}^{[\mu]}$ 
  satisfies condition ({\bf a}).
  For the lattice 
  $
  \Lambda_{\pmb{\alpha}}^{[\mu]}$ 
  the points
   \begin{equation}\label{pm1}
  \pmb{z}_\nu  =   \pm\mathcal{G}\pmb{m}_{\mu -\nu+1},\,\,\,\, \nu = 1,...,\mu
  \end{equation}
  are the best approximation points in the sense of this section, and the condition ({\bf b}) is satisfied.  We  choose  the signs $\pm$  in (\ref{pm1}) to have
  $
 z_{0,\nu} =+1,
$
and 
$
0<  z_{0,1} <  z_{0,2}<....< z_{ 0,\mu}.
$ We see that
\begin{equation}\label{invers}
z_{0,\nu} = T^{-d} L_{\mu-\nu+1},\,\,\,\,\,\,
|\underline{\pmb{z}}_\nu|_\infty = {M_{\mu-\nu+1}}{T}.
\end{equation}
If we take $ T \le M_\mu^{-1}$ we see that (\ref{0003}) is satisfied for all $ \nu = 1,...,\mu$.  We can take $ T$ small  enough to get
$$
|\underline{\pmb{z}}_\nu|_\infty z_{0,\nu+1}>
|\underline{\pmb{z}}_\nu|_\infty z_{0,\nu}=
T^{1-d} L_{\mu-\nu+1} M_{\mu-\nu+1}> 1\,\,\,\, \forall \nu = 1,...,\mu.
$$
 So the conditions of Lemma 2 are satisfied.
 Now we apply Lemma 2 and Lemma 1$'$  to show that 
 $$
 |\underline{\pmb{z}}_1|_\infty z_{0,2}\ge
 \frac{\Delta_2}{2\sqrt{2}\, d} =
  \frac{1}{2\sqrt{2}\, d} 
  \cdot
  \prod_{j=2}^d
   \frac{\Delta_j}{\Delta_{j+1}}  {\rm det}\,   \Lambda_{\pmb{\alpha}}^{[\mu]} =  
    \frac{1}{2\sqrt{2}\, d} 
  \cdot
  \prod_{j=2}^d
   \frac{\Delta_j}{\Delta_{j+1}} \gg_d
    \prod_{j=2}^d\left(
    \,\frac{z_{0,\nu_{j+1}-1}}{ 
    z_{0,\nu_{j+1}} } 
    \right)
      \cdot
         \prod_{j=2}^d\frac{1}{ 
    |\underline{\pmb{z}}_{\nu_{j+1}-1}|_\infty}.
 $$
 
 We have assumed ({\bf iii}), so
 $$
 \frac{L_{j+1}}{ 
    L_j } > \gamma >0 ,\,\,\,\,\, \forall j
    $$
    and by the first formula from (\ref{invers}) we get
 $$
 \frac{z_{0,\nu_{j+1}-1}}{ 
    z_{0,\nu_{j+1}} } =
    \frac{L_{\mu-\nu_{j+1}+2}}{ 
    L_{\mu -\nu_{j+1}+1} } > \gamma.
    $$
    From the other  hand,  for $j\ge 2$ we have $|\underline{\pmb{z}}_{\nu_{j+1}-1}|_\infty < |\underline{\pmb{z}}_{1}|_\infty$,
    because of $ \nu_{j+1}-1\ge \nu_3-1\ge \nu_2>\nu_1=1$ and
    $$
    |\underline{\pmb{z}}_{\nu_{j+1}-1}|_\infty = M_{\mu-\nu_{j+1}+2}T,\,\,\,\,\,
    |\underline{\pmb{z}}_1|_\infty = M_{\mu}T
    $$
    (see (\ref{invers})) and  $M_{\mu+1}< M_{\mu+\nu_{j+1}} $  (see (\ref{mioX})).
    We conclude with
     $$
 |\underline{\pmb{z}}_1|_\infty z_{0,2}\gg_{d,\gamma} 
 \frac{1}{|\underline{\pmb{z}}_{1}|_\infty^{d-1}},
 $$
 or 
 $$|\underline{\pmb{z}}_{1}|_\infty^d   z_{0,1}>
 \gamma
|\underline{\pmb{z}}_{1}|_\infty^d   z_{0,2}\gg_{d,\gamma}  1,
$$
as $  z_{0,1}/ z_{0,2}> \gamma$.
We apply (\ref{invers}) again to see that $ L_\mu M_\mu^d \gg_{d,\gamma}1$. The last inequality holds for all $\mu$  large enough and this means that $\pmb{\alpha}$ is badly approximable. 

\vskip+0.3cm

Now we suppose that 
 there exists $\nu_0$ and a proper linear subspace $ \mathcal{L}\subset \mathbb{R}^{d+1}$ 
  of dimension $3\le l={\rm dim}\, \mathcal{L} < d+1$ such that $ \pmb{m}_\nu \in \mathcal{L}$ for all  $\nu\ge \nu_0$.
  We 
  may suppose that  $\mathcal{L}$ 
  has the minimal dimension among all such subspaces.
 Then $\mathcal{L} $ is a rational subspace and
 inside $\mathcal{L} $  we have an irrational subspace
  $$
  \mathcal{L}_1=
  \{ \pmb{x}= (x_0,x_1,...,x_d) \in \mathcal{L}:\,\, x_0+ x_1\alpha_1+...+x_d\alpha_d = 0\}\subset \mathcal{L}.
  $$
    But then all the best approximations vectors $\pmb{m}_\nu$ will be all the best  approximation vectors of 
  the lattice $\mathcal{L}\cap \mathbb{Z}^{d+1}$ to $\mathcal{L}_1$ in the induced norm, and this means that  the values 
  $L_\nu$ are proportional to the values $ \rho (\pmb{m}_\nu, \mathcal{L}_1)$. 
  From the other hand
 the argument behind shows that the $(l-1)$-dimensional subspace  $  \mathcal{L}_1$
  is badly approximable in $\mathcal {L}$, that is
  $$
  \inf_{\pmb{m}\in \mathcal{L}\cap \mathbb{Z}^{d+1} \setminus \{ \pmb{0}\}}
  \rho (\pmb{m}, \mathcal{L}_1)    |\pmb{m}|^{l-1}>0.
  $$
  But then all the best approximations vectors $\pmb{m}_\nu$ will be all the best  approximation vectors of 
  the lattice $\mathcal{L}\cap \mathbb{Z}^{d+1}$ to $\mathcal{L}_1$ in the induced norm, and this means that  the values 
  $L_\nu$ are proportional to the values $ \rho (\pmb{m}_\nu, \mathcal{L}_1)$ and hence
   $$
  \inf_{\nu}  
 L_\nu  |\pmb{m}_\nu|^{l-1}
  >0
 . $$
This is not possible, because for $ l \le d$ this contradicts  (\ref{bo1}).

\vskip+0.3cm
So the proof is completed.$\Box$

\vskip+0.3cm

{\bf Remark 3.}\, {\it
In the last part of the proof we deal with the situation when the subspace of best approximations  for a linear form 
has dimension smaller than $d+1$.  In particular we proved that this is not possible  for badly approximable $\pmb{\alpha}$.
 Such type of problems    were discussed in a recent paper \cite{Schl}.}
 
\vskip+0.3cm

   \vskip+0.3cm
   {\bf  9. Construction of approximations in two-dimensional subspace.}
   
      \vskip+0.3cm
   
     The following obvious lemma will be very useful.
   \vskip+0.3cm

   {\bf Lemma 3.}\,{\it 
   Let  $ \pmb{v} = (p, b_{1},b_2)\in\mathbb{Z}^3$, $p\ge 1$  be a primitive integer vector
   and $\pmb{V} = 
   \left( \frac{b_1}{p},\frac{b_2}{p}\right)$ be the corresponding rational vector.
  Suppose that 
   $ \delta = \delta (  \pmb{v} ) =\frac{1}{2p^2}$.
   Then
 for all  
   $\pmb{\rm  x} $ under the condition
   $$
    |\pmb{\rm  x} -   \pmb{V}|_\infty
  <\delta
  $$
  the
  vector  $ \pmb{v} $ is a best approximation vector for $\pmb{\rm  x} $.
}
\vskip+0.3cm

Proof.  Let us assume for two independent vectors
$
 \pmb{v} = (p, b_{1},b_2)
 $
 and
 $
 \pmb{v}' =(p',b_1',b_2')\in\mathbb{Z}^3
 $
 with 
 $0<p'\le p$ the induced vectors
  $\pmb{V} = 
   \left( \frac{b_1}{p},\frac{b_2}{p}\right)$
   and
     $\pmb{V}' = 
   \left( \frac{b_1'}{p'},\frac{b_2'}{p'}\right)$
    both have  distance smaller that  $\delta$ from $\pmb{x}$.
    Then
    $$
    |  \pmb{V}  -   \pmb{V}' |_\infty
 \le
     |  \pmb{V}  -   \pmb{x} |_\infty 
     +
       |  \pmb{V}'  -   \pmb{x} |_\infty <\frac{1}{p^2}.
       $$
       On the other hand, since by linear independence
$ \pmb{V}  \neq  \pmb{V}' $       
   and both coordinates in the difference $ \pmb{V}  -   \pmb{V}'$  have common denominator $pp'\le p^2$,
   we have the reverse bound
   $ |  \pmb{V}  -   \pmb{V}' |_\infty\ge p^{-2}$, and this is a contradiction. Hence 
   $\pmb{v}'$ is linearly dependent to $\pmb{v}$. Finally since $\pmb{v}$ is primitive, there is no such integer vector 
   $\pmb{v}'\neq \pmb{v}$ with $ p'\le p$.
   $\Box$
 
   \vskip+0.3cm

  \vskip+0.3cm
   
   {\bf Lemma 4.}\, {\it  
   Suppose that two independent integer points 
   $$\pmb{v}_0 = (p_0, b_{1,0}, b_{2,0}),\,\,\,\pmb{v}_1 =(p_1, b_{1,1}, b_{2,1})\in \mathbb{Z}^3$$
    with 
   \begin{equation}\label{poo}
   p_1>p_0 \ge 1 
   \end{equation}
   and the corresponding rational points
    $$
    \pmb{V}_0 = \left(\frac{b_{1,0}}{p_0}, \frac{b_{2,0}}{p_0}\right),\,\,\,\pmb{V}_1 =\left(\frac{b_{1,1}}{p_1}, \frac{b_{2,1}}{p_1}\right)\in \mathbb{Q}^2
    \cap [0,1]^2
    $$
   satisfy the  following conditions.


     \vskip+0.3cm
   \noindent
   {\rm ({\bf i})} the lattice
   $\Lambda = \langle \pmb{v}_0,\pmb{v}_1\rangle_\mathbb{Z}$ is complete, that is 
   $$
   \langle \pmb{v}_0,\pmb{v}_1\rangle_\mathbb{Z} =
    \pi \cap \mathbb{Z}^3
$${
    where }
    $$
    \pi =  \langle \pmb{v}_0,\pmb{v}_1\rangle_\mathbb{R} 
    $$ 
    is a two-dimensional plane spanned by 
     $ \pmb{v}_0$ and $\pmb{v}_1$;
     by
     $ \Delta $ we denote
  the fundamental volume of two-dimensional lattice   
   $\Lambda = \langle \pmb{v}_0,\pmb{v}_1\rangle_\mathbb{Z}$;
   
     \vskip+0.3cm
   \noindent
   {\rm ({\bf ii})}    points $\pmb{V}_0$ and $\pmb{V}_1$ satisfy 
   \begin{equation}\label{x1}
    |\pmb{V}_0 -\pmb{V}_1|_\infty
   \le \frac{1}{2} \, \min \left( \frac{1}{p_0\Delta}, \delta  (\pmb{v}_0)\right),
   \end{equation}  
   where $\delta  (\pmb{v}_0)$ is defined in Lemma 3.
   \vskip+0.3cm

  Consider the vectors  $\pmb{v}_i =(p_i, b_{1,i}, b_{2,i}),   2\le i \le k$  defined recursively by
       \begin{equation}\label{x2}
  \pmb{v}_{i} =  \pmb{v}_{i-1}+ \pmb{v}_{i-2}, 
   \end{equation}  
    and the corresponding rational points 
 $$ \pmb{V}_i =  \left(\frac{b_{1,i}}{p_i}, \frac{b_{2,i}}{p_i}\right)\in \mathbb{Q}^2
   $$
    such that  
       \begin{equation}\label{u}
 p_k \ge   \varkappa = \varkappa ( \pmb{v}_0, \pmb{v}_1)=
 \max
 \left( \Delta^2,
 \sqrt{
 \frac{\Delta}{\delta(\pmb{v}_0)}
 }
 ,
 \sqrt{
 \frac{p_1\Delta}{|p_0(\pmb{V}_2-\pmb{V}_1)|_\infty}
 \left(
 1+\frac{p_1}{p_0}\right)
 }
 \right)
   \end{equation} 

 Then    for any $\pmb{\rm  x} = (x_1,x_2) \in \mathbb{R}^2$ satisfying 
   \begin{equation}\label{x3}
 |\pmb{\rm  x} -   \pmb{V}_k|_\infty
  \le \frac{\Delta}{100p_k^2}
   \end{equation}   
   either
  \begin{equation}\label{x5}
    \pmb{v}_0,  \pmb{v}_1,...  , \pmb{v}_{k-2},\pmb{v}_{k-1}, \pmb{v}_k,
   \end{equation}
   or 
   \begin{equation}\label{x51}
     \pmb{v}_0,  \pmb{v}_1,... ,  \pmb{v}_{k-2}, \pmb{v}_k
   \end{equation}
   is the sequence of all consecutive best approximation vectors  from
   $\pmb{v}_0$ to  $\pmb{v}_k$, that is all the best approximation vectors $\pmb{z}  =(q,a_1,a_2) $ to $\pmb{\rm x}$ with
   $p_0\le q\le p_k$. 
   
   Moreover for every $\pmb{\rm  x} $ under the consideration we have
   \begin{equation}\label{laci}
   \frac{
   |p_i \, \pmb{\rm  x} -\underline{\pmb{v}}_i|_\infty
   }
   {
   |p_{i-1} \, \pmb{\rm  x} -\underline{\pmb{v}}_{i-1}|_\infty  
   }
   \ge
   \frac{1}{4}
   ,
   \,\,\,\,\,\,\,\,\,
   i= 1,2,..., k-1.
   \end{equation}

   }
\vskip+0.3cm 

Proof. 
Let us start with any $\pmb{\rm x}$  satisfying (\ref{x3}).
For  $0\le i \le k$ consider points
$$\pmb{Z}_i = (p_i, p_ix_1,p_ix_2)\,\,\,\,\,
\text{and}\,\,\,\,\,
\frak{z}_i =  \left(p_i, p_i \frac{b_{1,k}}{p_k},p_i\frac{b_{2,k}}{p_k}\right)
$$  
and the remainder vectors
$$
\pmb{\eta}_i =  \pmb{Z}_i -\pmb{v}_i\,\,\,\,\,
\text{and}
\,\,\,\,\,
\frak{y}_i =  \frak{z}_i -  \pmb{v}_i.
$$
More generally, for a vector $ \pmb{v} = (p,b_1,b_2) \in    \pi
$
we write
$$
\frak{y} (\pmb{v} ) = 
\left(0, p \frac{b_{1,k}}{p_k}- b_1,p\frac{b_{2,k}}{p_k} - b_2\right).
$$
We should note here that
as all the vectors $\frak{y} (\pmb{v} )$ are parallel, their sup-norms
$ |\underline{\frak{y}} (\pmb{v} )|_\infty$
 are proportional to Euclidean norms $ |\frak{y} (\pmb{v} )|$, that is 
 for vectors  $ \pmb{v} = (p,b_1,b_2) ,\pmb{v}' = (p',b_1',b_2') \in    \pi$ we have
 \begin{equation}\label{o1}
 \frac{ |\underline{\frak{y}} (\pmb{v} )|_\infty}{
 |\underline{\frak{y}} (\pmb{v}' )|_\infty} =
 \frac{
  |\frak{y} (\pmb{v} )|}{
   |\frak{y} (\pmb{v}' )|}.
 \end{equation}
 It is clear that 
  $$
  |\underline{\frak{y}} (\pmb{v} )|_\infty
 \ge
  |\frak{y} (\pmb{v} )| /\sqrt{2}.
 $$

From (\ref{x2}) it follows that
$$
\pmb{\eta}_{i+1} = \pmb{\eta}_{i}+\pmb{\eta}_{i-1}\,\,\,\,\,
\text{and}
\,\,\,\,\,
\frak{y}_{i+1} =  \frak{y}_i + \frak{y}_{i-1}.
$$
In addition we may note that vectors 
$\frak{y}_i$ are parallel and
$$
\frak{y}_{i} = - \frac{|\frak{y}_{i}|}{|\frak{y}_{i-1}|} \cdot \frak{y}_{i-1}.
$$
So 
$$
  |\frak{y}_{k}|=0,\,\,\,\,\,
|\frak{y}_{k-1}|=|\frak{y}_{k-2}|,\,\,\,\,\,
|\frak{y}_{i-1}| = |\frak{y}_{i}|+|\frak{y}_{i+1}|,
$$
and we can write the ratio 
$
\frac{|\frak{y}_{i-1}|}{|\frak{y}_{i}|}  
$
as the continued fraction and get the estimates
$$
   \frac{|\frak{y}_{i-1}|}{|\frak{y}_{i}|} = [\, \underbrace{1;1,...,1}_{k-i}\,] \le 2,\,\,\, \,\,\,\,\,1\le i\le k-1
$$
and
$$
   \frac{|\frak{y}_{i-1}|}{|\frak{y}_{i}|} \ge \frac{3}{2} ,\,\,\,\,\,\,\, 1\le i\le k-2
$$
So
  by (\ref{o1}),
\begin{equation}\label{x6}
 \frac{ |\underline{\frak{y}}_i|_\infty}{
 |\underline{\frak{y}}_{i-1}|_\infty} =
 \frac{
  |\frak{y}_i|}{
   |\frak{y}_{i-1}|}
     \ge \frac{1}{2}   ,\,\,\,\,\,\,\,\,\,\, 1\le i \le k-1,
   \end{equation}
   and
   \begin{equation}\label{x61}
 \frac{ |\underline{\frak{y}}_i|_\infty}{
 |\underline{\frak{y}}_{i-1}|_\infty} 
  \le \frac{2}{3}   ,\,\,\,\,\,\,\,\,\,\, 1\le i \le k-2.
   \end{equation}
   We should note that the point $\pmb{V}_k$ belongs to the segment with endpoints $\pmb{V}_0$, $\pmb{V}_1$ which belong to   the plane $\pi$.
  From (\ref{x1})  
  we see that  $|\pmb{V}_k- \pmb{V}_0|< \frac{\delta (\pmb{v}_0)}{2}$. So by Lemma 3, $\pmb{v}_0$ is a best approximation vector to $\pmb{V}_k$.
  Moreover, 
  the integer lattice $\mathbb{Z}^3$ splits into two-dimensional sublattices parallel to $\pi$. The Euclidean distances between  the corresponding neighboring  two-dimensional planes is 
  equal to $\Delta^{-1}$. So from (\ref{x1})  we see that for any integer point $(p'',b_1'',b_2'')\in \mathbb{Z}^3 \setminus \pi$
  one has 
  $|p''\pmb{V}_k - \underline{\pmb{b}}''|_\infty > |p_0 \pmb{V}_k - \underline{\pmb{b}}|_\infty$.
  So
  we deduce that all the best approximations to   $\pmb{V}_k$
 with denominator greater than $p_0$ lie in the plane $\pi$. 
 
 Now we consider  an approximation to $\pmb{V}_k$ from subspace $\pi$.
 For any $ i = 1,..., k$ 
 the points $\pmb{v}_{i-1}, \pmb{v}_i\in \pi $ form a basis of $\Lambda$. Moreover  the points 
 $\pmb{v}_{i-1}, \pmb{v}_{i}$ lie on the opposite sides from the line  $\langle \pmb{v}_k\rangle_\mathbb{R}$.
(Here we should note that for the case $ i = k$ the point $\pmb{v}_{i}$ lies  just on the line 
$\langle \pmb{v}_k\rangle_\mathbb{R}$, however our argument remains valid.)
 So there is no vectors $ \pmb{v} = (p,b_1,b_2) \in \pi$ satisfying
 $$
 p_{i-1} < p< p_{i},\,\,\,\,\,
 \text{and}
 \,\,\,\,\,
 \pmb{v} = \lambda \pmb{v}_{i-1}+\mu \pmb{v}_i,\,\,\,
 \lambda \in \{0,1\} ,\,\, \mu \in \mathbb{Z}.
 $$
 We see that for any vector  $ \pmb{v} = (p,b_1,b_2) \in \pi$ 
 with  $
 p_{i-1} < p< p_{i}
 $ we have
 $$
  \pmb{v} = \lambda \pmb{v}_{i-1}+\mu \pmb{v}_i,\,\,\,
 \lambda \neq 0,1 ,\,\, \mu \in \mathbb{Z}.
 $$
 Consider the lines
  $$
 \ell =\ell(\lambda) = \{ \pmb{x} = \lambda \pmb{v}_{i-1} + \mu \pmb{v}_i, \mu \in \mathbb{R}\}\,\,\,\,\,
 \lambda\in \mathbb{Z}.
 $$
 We should note that if  points $ \pmb{v} = (p,b_1,b_2)
 \in \ell( \lambda) 
 $ and $\pmb{v}' = (p, b_1',b_2')  \in \ell( \lambda') $ with the same first coordinate $p \in (p_{i-1}, p_i)$  belong to two parallel  lines
 $
 \ell(\lambda) 
 $
 and $    \ell(\lambda') 
  $ with  integers $ \lambda \neq \lambda'$
 then
 \begin{equation}\label{oi}
 |{\frak{y}}  ( \pmb{v}  -\pmb{v}') | =
 |\pmb{v}  -\pmb{v}'| 
 \ge
 \min_{ \pmb{v} \in \ell (0) , \pmb{v}'\in \ell (1)}  |\pmb{v}  -\pmb{v}' |=
    \sigma_i |\frak{y}_{i-1}| 
 ,\,\,\,\,
  \text{where}\,\,\,\,
    \sigma_i  =
   \left(1+
  \frac{
  |\frak{y}_{i}|}{  
 |\frak{y}_{i-1}|} \cdot\frac{p_{i-1}}{p_i}\right).
 \end{equation}
 We would like to give a comment on the last equality in (\ref{oi}). To obtain this inequality one should note that 
 $$
  \min_{ \pmb{v} \in \ell (0) , \pmb{v}'\in \ell (1)}  | \pmb{v}  -\pmb{v}' |=
   |\pmb{v} _{i-1}    -\frak{z}_{i-1}' |=
   |\pmb{v} _{i-1}  -\frak{z}_{i-1}|
  +
   |\frak{z}_{i-1}  -\frak{z}_{i-1}' |,
   $$
   where
   $$
   \frak{z}_{i-1}' = \left( p_{i-1}, p_{i-1} \frac{b_{1,i}}{p_i}, 
    p_{i-1} \frac{b_{2,i}}{p_i}
 \right).
 $$
 But   $ |\pmb{v} _{i-1}  -\frak{z}_{i-1}| = |\frak{y}_{i-1}|
 $
 and 
 $ 
  |\frak{z}_{i-1}  -\frak{z}_{i-1}' |=
  \frac{p_{i-1}}{p_i} |\frak{y}_{i}|$,
  and (\ref{oi}) follows. 
  
  So (\ref{oi}) shows that for all  $i = 1,..., k$  and for all $\pmb{v} = (p,b_1,b_2) \in \pi $ with 
  $ p_{i-1} <p<p_i$ we have
  \begin{equation}\label{oi1}
 |\frak{y}(\pmb{v})|
 \ge \sigma_i |\frak{y}_{i-1}|.
 \end{equation}
 Now from  (\ref{x6}) and 
 the inequality $ \frac{p_{i-1}}{p_i} \ge \frac{1}{2},  2\le i \le k-1$ we see that 
 $$ 
 \sigma_1 \ge 1+ \frac{p_0}{2p_1} > 1:\,\,\,\,\,
 \sigma_i \ge \frac{5}{4} 
 \,\,\,\,\,
 \text{for}
 \,\,\,\,\,
 2\le i \le k-1
 $$
 and (\ref{oi1}) transforms into
 \begin{equation}\label{y10}
{ |\underline{\frak{y}} (\pmb{v} )|_\infty}
   \ge
  \left(
   1+ \frac{p_0}{2p_1}
   \right)
   |\underline{\frak{y}}_{0}|_\infty,
   \,\,\,\,\,
   \text{for all}\,\,\,\,\,
 \pmb{v} = (p,b_1,b_2) \in \pi\,\,\,\,\,
 \text{
 with }\,\,\,\,\,
 p_{0} < p< p_{1};
 \end{equation}

 \begin{equation}\label{y1}
{ |\underline{\frak{y}} (\pmb{v} )|_\infty}
   \ge
   \frac{5}{4} |\underline{\frak{y}}_{i-1}|_\infty
   \,\,\,\,\,
   \text{for all}\,\,\,\,\,
 \pmb{v} = (p,b_1,b_2) \in \pi\,\,\,\,\,
 \text{
 with }\,\,\,\,\,
 p_{i-1} < p< p_{i}, \,\,\,\,\, 2\le i \le k-1.
 \end{equation}
By the same argument
 \begin{equation}\label{y11}
{ |\underline{\frak{y}} (\pmb{v} )|_\infty}
   \ge
       |\underline{\frak{y}}_{k-2}|_\infty+    |\underline{\frak{y}}_{k-1}|_\infty
   \frac{p_{k-2}}{p_{k-1}}\ge
   \frac{3}{2} |\underline{\frak{y}}_{k-2}|_\infty
   \,\,
   \text{for all}\,\,
 \pmb{v} = (p,b_1,b_2) \in \pi\,\,
 \text{
 with }\,\,\,\,\,
 p_{k-1} < p< p_{k}.
 \end{equation}


Now we  see  that  (\ref{x51}) 
 is the sequence of all best approximation vectors to  $\pmb{V}_k$ with denominators between  $p_0$ and $p_k$.
 As for the point $\pmb{v}_{k-1}$, it is not a best approximation vector because
 $|\frak{y}_{k-1}|=|\frak{y}_{k-2}|$ and so 
 $ { |\underline{\frak{y}}_{k-1}|_\infty}={
 |\underline{\frak{y}}_{k-2}|_\infty}$
 and $ p_{k-1} > p_{k-2}$.

 Now we need to estimate 
 $|\underline{\frak{y}}_{i}|_\infty , i = 0,..., k-1$ from below.
We consider   the lattice $\Lambda$  and the parallelogram
 $$
 \Pi = 
 \left\{(x,y_1,y_2) \in \mathbb{R}^3:
 \,\,
 |x|\le p_{i+1},\,\,
 \max_{j=1,2} 
    \left|
   x\frac{b_{j,k}}{p_k} - y_{j}  
   \right|\le
  |\underline{\frak{y}}_{i}|_\infty
   \right\} \cap \pi.
 $$
 For its area we have
 $$
 \Delta \le {\text{area} (\Pi ) \le 
 2\sqrt{2}   |\underline{\frak{y}}_{i}|_\infty 
 \times \sqrt{1+ \left( \frac{b_{1,k}}{p_k}\right)^2+ \left( \frac{b_{2,k}}{p_k}\right)^2}} \,\,p_{i+1}
 \le2\sqrt{6} \,  |\underline{\frak{y}}_{i}|_\infty\,  p_{i+1}.
 $$ 
So
  \begin{equation}\label{l2}
 |\underline{\frak{y}}_{i}|_\infty \ge 
    \frac{\Delta}{ 2\sqrt{6} p_{i+1}},\,\,\,\,\,
     i = 0,..., k-1.
 \end{equation}

 Now we prove the statement of the lemma  about points $\pmb{x}$ under the condition (\ref{x3}).
 From (\ref{x3}) we see that  
\begin{equation}\label{x10w}
 |p\,\pmb{\rm  x} -   p\pmb{V}_k|_\infty=
  \max_{j=1,2} \left|
  px_j - p\frac{b_{j,k}}{p_k}
  \right| \le  \frac{\Delta}{100p_k} \cdot \frac{p}{p_k}
 ,
\end{equation}
and in particular  for $p \le p_k$ one has
\begin{equation}\label{x10}
 |p\,\pmb{\rm  x} -   p\pmb{V}_k|_\infty
 \le  \frac{\Delta}{100p_{k}} .
\end{equation}

First of all we show
that 
the vectors $\pmb{v} = (p,b_1,b_2)  \in \mathbb{Z}^3 \setminus \Lambda$  with $ p >p_0$ cannot be best approximation vectors for $\pmb{\rm x}$. Indeed, 
  $ \pmb{V}_k $ belongs to the segment with endpoints   $ \pmb{V}_0, \pmb{V}_1 $ and 
  inequality  (\ref{x1}) 
show that
$$
|  \pmb{V}_k -  \pmb{V}_0| < \frac{1}{2p_0\Delta}.
$$
Now (\ref{x10}) together with the inequality (\ref{u}) written as
$p_k \ge  \Delta^2$
give the bound
$$
 |p_0\,\pmb{\rm  x} -   \underline{\pmb{v}}_0|_\infty
 \le
  |p_0\,\pmb{V}_k -   \underline{\pmb{v}}_0|_\infty+
   |p_0\,\pmb{\rm  x} -  p_0\,\pmb{V}_k |_\infty
    \le  \frac{1}{2\Delta} +
   \frac{\Delta}{100p_{k}}
   \le   \frac{1}{\Delta}  \le
 |p\,\pmb{\rm  x} -   \underline{\pmb{v}}|_\infty.
$$

Then we show
that $\pmb{v}_0$ is a best approximation for  $\pmb{\rm  x}$.
Indeed,
from (\ref{x1}, \ref{x3}) and  (\ref{u}) in  the form
$ p_k \ge  \sqrt{\frac{\Delta}{\delta(\pmb{v}_0)}}$
we have
$$
|\pmb{V}_0 - \pmb{\rm  x} | \le
|\pmb{V}_k - \pmb{V}_0|+
|\pmb{V}_k - \pmb{\rm  x} | 
 \le {\delta{(\pmb{v}_0})}.
$$
So $\pmb{v}_0$ is the best approximation vector for $\pmb{\rm  x} $.

 Now we study approximation to  $  \pmb{\rm  x}$  by vectors $ \pmb{v}_i,  i = 0,1,...,k$.

From the  triangle inequality  and (\ref{x10})  for vectors
$  p_i\, \pmb{\rm  x}$ and  $  \underline{\pmb{v}}_i = (b_{1,i},b_{2,i})$
we deduce for
$
  |\underline{\frak{y}}_{i}|_\infty  =|p_i \pmb{V}_k -   \underline{\pmb{v}}_i
|_\infty$ the inequalities
\begin{equation}\label{pppp}
  |\underline{\frak{y}}_{i}|_\infty 
 -  \frac{\Delta}{100p_{k}}
\le
 |p_i\,\pmb{\rm  x} -   \underline{\pmb{v}}_i|_\infty
  \le
  |\underline{\frak{y}}_{i}|_\infty +  \frac{\Delta}{100p_{k}}
  ,\,\,\,\,\,
  i = 0,..., k
.
\end{equation}
We should note that from (\ref{l2}) we have
\begin{equation}\label{ppp}
\frac{\Delta}{100 p_k  |\underline{\frak{y}}_{i}|_\infty } \le \frac{\sqrt{6}}{50} ,\,\,\,\,\,\, 0\le i \le k-1.
\end{equation}
So
from the last inequality and (\ref{x61})    for $ i =1,...,k-2$ we get
\begin{equation}\label{fu}
\frac{
 |p_i\,\pmb{\rm  x} -   \underline{\pmb{v}}_i|_\infty}{
  |p_{i-1}\,\pmb{\rm  x} -   \underline{\pmb{v}}_{i-1}|_\infty} \le
 \frac{
  |\underline{\frak{y}}_{i}|_\infty 
  }
  {
  |\underline{\frak{y}}_{i-1}|_\infty 
  } \cdot
   \frac{1+
  \frac{\Delta}{100p_{k} |\underline{\frak{y}}_{i}|_\infty }} 
  {1-
  \frac{\Delta}{100p_{k}{  |\underline{\frak{y}}_{i-1}|_\infty }} 
  }
 \le
  \frac{3}{4}
\,\,\,\,\, \,\,\,\,\,
1\le i \le k-2.
\end{equation} 
In addition from  (\ref{x3}), (\ref{pppp}) and (\ref{ppp}) we deduce
\begin{equation}\label{fu1}
\frac{
 |p_k\,\pmb{\rm  x} -   \underline{\pmb{v}}_k|_\infty}{
  |p_{\nu}\,\pmb{\rm  x} -   \underline{\pmb{v}}_{\nu}|_\infty} \le
  \frac{\Delta}{100p_k |\underline{\frak{y}}_{\nu}|_\infty }\cdot
  \frac{1}{1-   \frac{\Delta}{100p_k |\underline{\frak{y}}_{\nu}|_\infty }}<\frac{1}{2},\,\,\,\,\,\,\,
  \nu = k-2, k-1.
\end{equation}

Let us show that there is no best approximations $ \pmb{v}= (p, b_1,b_2)$ with
$p_{i-1} < p <p_{i}$ for all $ i = 1,2,...,k$.

First of all we consider  the case $ i = 1$ that is $p_0 < p<p_1$.
In this case we will take into account the inequality
$$
|\underline{\frak{y}}_0|_\infty \ge |p_0 (\pmb{V}_2-\pmb{V}_0)|_\infty,
$$
as well as the inequalities
$$
\left|
|
p_0 
\,\pmb{\rm  x} -
\underline{\pmb{v}}_0
|_\infty -
|\underline{\frak{y}}_0
|
\right|_\infty
\le \frac{p_0\Delta}{100 p_k^2},\,\,\,\,\,\,\,\,\,
\left|
|
p
\,\pmb{\rm  x} -
\underline{\pmb{v}}
|_\infty -
|\underline{\frak{y}} (\pmb{v})
|_\infty
\right|
\le \frac{p_1\Delta}{100 p_k^2}
,
$$
which follow from (\ref{x10w}). Three last inequalities together with (\ref{u}) in the form
$$ p_k \ge \varkappa\ge \sqrt{
\frac{p_1\Delta}{100|p_0(\pmb{V}_2-\pmb{V}_1)|_\infty} \left( 1+\frac{p_1}{p_0}\right)}
$$
and (\ref{y10})
lead to
$$
\frac{|p
\,\pmb{\rm  x} -\underline{\pmb{v}}
|_\infty }{
|p_0
\,\pmb{\rm  x} -\underline{\pmb{v}}_0
|_\infty }\ge
\frac{|
\underline{\frak{y}}(\pmb{v})|_\infty - \frac{p_1\Delta}{100 p_k^2}}
{
|
\underline{\frak{y}}_0|_\infty +\frac{p_0\Delta}{100 p_k^2}
}\ge
\frac{
|
\underline{\frak{y}}(\pmb{v})|_\infty }
{
|
\underline{\frak{y}}_0|_\infty  
}
\cdot
\frac{ 1 - \frac{p_1\Delta}{100p_k^2|
\underline{\frak{y}}_0|_\infty}}
{
1 +\frac{p_0\Delta}{100p_k^2 |
\underline{\frak{y}}_0|_\infty}
}>1,
$$
and we proved everything  what we need in the case $ p_0 <p<p_1$.

Next, suppose that 
$p_{i-1} <p<p_{i}\,\,\,\,\,\text{and} \,\,\,\,\, 2\le  i \le k-1.$
 Then (\ref{y1}) and (\ref{ppp})  give
 $$
\frac{|p
\,\pmb{\rm  x} -\underline{\pmb{v}}
|_\infty }{
|p_{i-1}
\,\pmb{\rm  x} -\underline{\pmb{v}}_{i-1}
|_\infty }\ge 
\frac{
|
\underline{\frak{y}}(\pmb{v})|_\infty }
{
|
\underline{\frak{y}}_{i-1}|_\infty  
}
\cdot
\frac{ 1 - \frac{p_1\Delta}{100p_k^2|
\underline{\frak{y}}_{i-1}|_\infty}}
{
1 +\frac{p_1\Delta}{100 p_k^2 |
\underline{\frak{y}}_{i-1}|_\infty}
}>1,
$$
and  everything is done in  the case $ p_1<p<p_{k-1}, p \neq p_i$ also.

By similar argument using  (\ref{y11}) and (\ref{ppp}) for $p_{k-1}<p <p_k$ we see that
 $$
\frac{|p
\,\pmb{\rm  x} -\underline{\pmb{v}}
|_\infty }{
|p_{i-2}
\,\pmb{\rm  x} -\underline{\pmb{v}}_{i-1}
|_\infty }>  1.
$$
We see from (\ref{fu},\ref{fu1}) and the lower bounds for 
$|p
\,\pmb{\rm  x} -\underline{\pmb{v}}
|_\infty$ that 
 $\pmb{v}_0,\pmb{v}_1,...,\pmb{v}_{k-2}$ and $\pmb{v}_k$  are the best approximation vectors for $\pmb{\rm  x}$, and $\pmb{v}_{k-1}$ may be a best approximation vector or may be not.
So all the best approximations for $\pmb{\rm x}$ form either the sequence (\ref{x5})
or the sequence (\ref{x51}).
 
 \vskip+0.3cm
 To finish the proof of Lemma 4 we need to show (\ref{laci}). But this can be done analogously to (\ref{fu}), as from (\ref{x6}) and (\ref{ppp}) we see that 
 $$
 \frac{
 |p_i\,\pmb{\rm  x} -   \underline{\pmb{v}}_i|_\infty}{
  |p_{i-1}\,\pmb{\rm  x} -   \underline{\pmb{v}}_{i-1}|_\infty} \ge
 \frac{
  |\underline{\frak{y}}_{i}|_\infty 
  }
  {
  |\underline{\frak{y}}_{i-1}|_\infty 
  } \cdot
   \frac{1-
  \frac{\Delta}{100p_{k} |\underline{\frak{y}}_{i}|_\infty }} 
  {1+
  \frac{\Delta}{100p_{k}{  |\underline{\frak{y}}_{i-1}|_\infty }} 
  }
\ge \frac{1}{4},
\,\,\,\,\, \,\,\,\,\,
1\le i \le k-1.
$$

  $\Box$

     \vskip+0.3cm

  {\bf 10. Three-dimensional subspaces.}
  
  \vskip+0.3cm
  {\bf Lemma 5.}\, {\it  
Consider two independent integer points
  $$
  \pmb{w}_0' =  (p_0', b_{1,0}', b_{2,0}'),
 \,\,\,\,\,
  \pmb{w}_0'' =  (p_0'', b_{1,0}'', b_{2,0}'')
  $$
 and the two-dimensional subspace $ \pi=\langle    \pmb{w}_0',   \pmb{w}_0''\rangle_\mathbb{R}$.
 Suppose that for the corresponding rational points we have
  $$
 \pmb{W}_0'=\left(\frac{b_{1,0}'}{p_0'},  \frac{b_{1,0}'}{p_0'}\right),\,\,\,\,\,
 \pmb{W}_0''=\left(\frac{b_{1,0}''}{p_0''},  \frac{b_{1,0}''}{p_0''}\right)\in [0,1]^2.
 $$
 Suppose that  
 $  \pmb{w}_0'$ and $  \pmb{w}_0''$
 form a basis of the lattice
 $\Lambda =\pi \cap \mathbb{Z}^3$, that is
 $$
 \Lambda =  \langle    \pmb{w}_0',   \pmb{w}_0''\rangle_\mathbb{Z}
 ,$$
 and 
 $\Delta$  is the two-dimensional fundamental volume of $\Lambda$.
     Suppose  
     that parameters $\gamma_1$ and $\gamma_2$ satisfy the inequalities
   \begin{equation}\label{gammaI}
   \gamma_2 \ge \gamma _1^ 2,\,\,\,\,\, \gamma_1 \ge 50.
   \end{equation}
 Consider the  point
 $$
 \pmb{w}_0 = 
  \pmb{w}_0'+\pmb{w}_0'' = (p_0, b_{1,0},b_{2,0})
  $$
  and the corresponding rational point
   and 
    $\pmb{W}_0= \left(\frac{b_{1,0}}{p_0},  \frac{b_{1,0}}{p_0}\right)\in [0,1]^2$.
Suppose that
 \begin{equation}\label{del}
 p_0 \ge \gamma_1 \Delta^2.
 \end{equation}
     
    Let $\pmb{n}$ be an orthogonal vector to $ \pi $ and $ |\pmb{n}| = 1$.
  Consider the point
  \begin{equation}\label{norma}
  \frak{x}_0 =(x_0,y_{1,0}, y_{2,0}) = \pmb{w}_0 + \pmb{n} \cdot \frac{\Delta}{\gamma_1 p_0}\in \mathbb{R}^3
  \end{equation}
  and the corresponding two-dimensional point
  \begin{equation}\label{norma1}
  \pmb{\rm x}_0 =  (x_{1,0},x_{2 ,0})=\left(
  \frac{y_{1,0}}{x_0}, \frac{y_{2,0}}{x_0}
  \right) \in \mathbb{R}^2.
\end{equation}
  Suppose that for all
  $
  \pmb{\rm x} = (x_1, x_2)\in\mathbb{R}^2$ satisfying
    \begin{equation}\label{o}
    |\pmb{\rm  x} -   \pmb{\rm  x}_0|_\infty=
  \max_{j=1,2}
  |x_j - x_{j,0}| \le \frac{\Delta}{\gamma_1 p_0^2}
  \end{equation}
the vector 
$
  \pmb{w}_0
  $ is a best approximation vector to $\pmb{\rm x}$.  
  
  Then there exists an integer point $ \pmb{w}_1 = (p_1,b_{1,1}, b_{2,1}) $ with the following properties:

  \vskip+0.3cm 
 \noindent 
{\rm ({\bf i})}  $\pmb{w}_1$ belongs to the affine subspace $\pi_1=  \pi + \frac{1}{\Delta}\cdot \pmb{n}$;

  \vskip+0.3cm

 \noindent 
 {\rm ({\bf ii)} } both triples 
 $$
  \pmb{w}_0',
 \,\,\,\,\,
  \pmb{w}_0,
  \,\,\,\,\,
  \pmb{w}_1
  $$
  and
   $$
  \pmb{w}_0'',
 \,\,\,\,\,
  \pmb{w}_0,
  \,\,\,\,\,
  \pmb{w}_1
  $$
  form bases in $\mathbb{Z}^3$;

    \vskip+0.3cm
   \noindent 
   {\rm ({\bf iii)} } 
   vectors 
    $
     \pmb{w}_0,
  \,
  \pmb{w}_1
  $
  form a basis of the two-dimensional lattice
  $$
  \Lambda_1 = \langle   \pmb{w}_0,
  \pmb{w}_1\rangle_\mathbb{R}\cap \mathbb{Z}^3
  $$
  with two-dimensional fundamental volume $\Delta_1$;
  
    \vskip+0.3cm
   \noindent 
   {\rm ({\bf iv)} }  the inequalities 
   \begin{equation}\label{w1}
 \left(
 {\gamma_1}-\frac{2}{\gamma_1}
 \right)\cdot\left( \frac{p_0}{\Delta}\right)^2\le
  p_1  \le \left(
 {\gamma_1}+\frac{2}{\gamma_1}
 \right) \cdot \left( \frac{p_0}{\Delta}\right)^2
   \end{equation}
   and
\begin{equation}\label{w2}
\frac{1}{4} \cdot \frac{p_0}{\Delta}
\le 
   \Delta_1 \le12 \cdot  \frac{p_0}{\Delta}
   \end{equation}
   are valid \footnote{It is important that  the constants  in 
   (\ref{w2}) do not depend on $\gamma_1$.};
   
   { 
     \vskip+0.3cm
     \noindent 
   {\rm ({\bf v)} }
   define
    $ \pmb{W}_1 = (\frac{b_{1,1}}{p_1}, \frac{b_{2,1}}{p_1}) $, then
    for any $\pmb{\rm x}= (x_1,x_2) \in \mathbb{R}^2$ satisfying
\begin{equation}\label{pat}
    |\pmb{\rm  x} -   \pmb{W}_1|_\infty=
  \max_{j=1,2}
  \left| x_j - \frac{b_{j,1}}{p_1} \right| \le \frac{\Delta}{ {\gamma_2 p_0 p_1}}
  \end{equation}  }
    either the vectors 
  $$
  \pmb{w}_0, \,\,\,
  \pmb{w}_1
  $$
  are two consecutive best approximation vectors to $\pmb{\rm x}$
  or the 
   vectors 
  $$
  \pmb{w}_0, \,\,\,\pmb{w}_1-\pmb{w}_0,\,\,\,
  \pmb{w}_1
  $$
  are three consecutive best approximation vectors to $\pmb{\rm x}$

 }
 
    \vskip+0.3cm
    {\bf Remark 4.}\, {\it For the point $\pmb{\rm x}_0  $ one has
    $$
   |\pmb{\rm  x}_0 -   \pmb{W}_0|
    \le
    \frac{\Delta}{\gamma_1 p_0^2}+
  (k-1)
    \sqrt{
1+ \left(\frac{b_{1,0}}{p_0}\right)^2 +
 \left(\frac{b_{2,0}}{p_0}\right)^2   + 
 \left(
\frac{\Delta}{\gamma_1 p_0^2}\right)^2
 }  \le
    \frac{2\Delta}{\gamma_1 p_0^2}
,$$
 where
 $ 1\le k =    \frac{p_0}{p_0 -  \frac{\Delta}{\gamma_1 p_0^2}
\sin \psi
}  $  
 and $ \psi $ is the angle between $\pmb{n} $ and $ \pmb{e} = (1,0,0)$.
 (We take into account that $ \frac{b_{j,0}}{p_0}\in [0,1]$.)
    }

    \vskip+0.3cm
    {\bf Remark 5.}\, {\it  From inequalities (\ref{gammaI},\ref{w1},\ref{w2}) it follows that
    $$
    \frac{1}{8}\cdot \frac{\Delta}{p_0}
    \le
    \frac{\Delta_1}{p_1}
    \le
    24 \cdot \frac{\Delta}{p_0}.
    $$
    }

    \vskip+0.3cm
    {\bf Remark 6.}\, {\it  From inequalities (\ref{w1}) and (\ref{dee}) it follows that
    $$
     p_1\ge \frac{\gamma_1}{2}\, p_0 \,\frac{p_0}{\Delta^2}\ge \frac{\gamma_1}{2}\, p_0.
    $$
    }
    
   \vskip+0.3cm
   Proof of Lemma 5. We should note that  the parallelogram
   $$
   \Pi = \{ \pmb{z}\in \mathbb{R}^3:\,\,\, \pmb{z}= \lambda\pmb{w}_0' + \mu\pmb{w}_0'',\,\,\, 0\le \lambda,\mu \le 1\}
   $$
   is a fundamental domain with respect to $\Lambda$ and the two-dimensional affine subspace  $\pi_1$ contains a lattice $\Lambda_1\subset \mathbb{Z}^3$ congruent to $\Lambda$.
   Then any shift of parallelogram $\Pi$ which belongs to $\pi_1$ contains an integer point. Consider the point
   $$
   \pmb{X}_0 =(X_0, Y_{1,0}, Y_{2,0})= \langle \pmb{x}_0\rangle_\mathbb{R} \cap \pi_1  
  $$
  and the parallelogram
  $
  \Pi+ \pmb{X}_0.
  $ By the discussion above it contains an integer point. We denote this point by
  $\pmb{w}_1 = 
 (p_1, b_{1,1}, b_{2,1})$. This is just the integer point what we need.
 Indeed, properties    {\rm ({\bf i)} }  and   {\rm ({\bf ii)} } 
 are clearly satisfied. As vector $\pmb{w}_0$ is primitive and there is no integer points between subspaces 
$\pi$ and $\pi_1$, property  {\rm ({\bf iii)} } is satisfied also. 
From the construction we see that 
 \begin{equation}\label{suomi}
X_0 = x_0 \cdot \frac{\gamma_1 p_0}{\Delta^2},\,\,\,\,\,
\,\,\,\,\,
|p_0 -x_0|\le \frac{\Delta}{\gamma_1p_0},
\,\,\,\,\,\,\,\,\,\,
|p_1-X_0|\le p_0.
\end{equation}
So
$$
\left|p_1-\frac{\gamma_1 p_0^2}{\Delta^2} \right|\le p_0 +\frac{1}{\Delta}< 2p_0
$$ 
and we get (\ref{w1}) by taking into account (\ref{del}).
To get (\ref{w2}) we
will estimate the area  $\Delta_1$ of  parallelogram
$$
\mathcal{P} = \{ \pmb{z} = a \pmb{X}_0+b\pmb{w}_1,\,\,\,\,\,
a,b\in [0,1)\}.
$$ 
 Observe that 
$$
\Delta_1=
{\rm area}\, 
\mathcal{P}
=
{\rm area}\, 
\mathcal{P}_0
+
\lambda_*\,
{\rm area}\, 
\mathcal{P}'
+
\mu_*\,
{\rm area}\, 
\mathcal{P}''
,
$$
with some $\lambda_*, \mu_* \in (-1,1)$,
where
$$
\mathcal{P}_0 = \{ \pmb{z} =a \pmb{w}_0+b \pmb{w}_1,\,\,\,\,\,
a,b\in [0,1)\}
$$
and
$$
\mathcal{P}' = \{ \pmb{z} = a\pmb{w}_0'+b \pmb{w}_1,\,\,\,\,\,
a,b \in [0,1)\},
\,\,\,
\mathcal{P}'' = \{ \pmb{z} = a \pmb{w}_0''+b \pmb{w}_1,\,\,\,\,\,
a,b\in [0,1)\}.
$$ 
It is clear that 
$$
{\rm area}\, 
\mathcal{P}'
, 
{\rm area}\, 
\mathcal{P}''
\le \Delta,
$$
and
$$
{\rm area}\, 
\mathcal{P}_0 = |\pmb{X}_0 | 
\rho (\frak{x}_0 , \langle\pmb{w}_0\rangle_\mathbb{R}) =
|\pmb{X}_0 | \cdot  \frac{\Delta}{\gamma_1 p_0} ,
$$
 {where}
 $$
\left(
p_0 - \frac{\Delta}{\gamma p_0}
\right)  \frac{\gamma p_0}{\Delta^2}
\le
|\pmb{X}_0 |
=|x_0 |\cdot \frac{\gamma_1p_0}{\Delta^2}
\le \left(\sqrt{3} p_0 +\frac{\Delta}{\gamma_1 p_0}\right) \frac{\gamma_1 p_0}{\Delta^2}
.
$$
So
$$
  |\pmb{X}_0 |\cdot \frac{\Delta}{\gamma_1 p_0}  - 2\Delta
\le
\Delta_1  \le
 |\pmb{X}_0 |\cdot
\frac{\Delta}{\gamma_1 p_0}  + 2\Delta  ,
$$
and together with  (\ref{w1}) the last two formulas   give  (\ref{w2}).

  To finish the proof it remains to  explain  {\rm ({\bf v})}. 
 
 If 
 $\pmb{\rm x}$  satisfies (\ref{pat})  then it     satisfies (\ref{o}).
 Indeed,   as $\pmb{w}_1 \in  \Pi+ \pmb{X}_0
 $,
 the point
 $\pmb{W}_1=\left(\frac{b_{1,1}}{p_1}, \frac{b_{2,1}}{p_1}
 \right)
 $
 belongs to a convex polygon with vertices
 $$
 \pmb{\rm x}_0 = (x_{1,0},x_{2,0}),
 \,\,\,\,\,\,\,\,\,\,\,\,\,\,\,\,\,\,\,\,\,\,\,\,\,\,\,\,\,\,\,\,\,\,\,\,\,\,\,\,\,\,\,\,\,\,
  \pmb{\rm x}_{0,0}=
  \left(\frac{Y_{1,0}+b_{1,0}}{X_0+ p_0}
 ,\frac{Y_{2,0}+b_{2,0}}{X_0+ p_0}
 \right),
 $$
 $$
  \pmb{\rm x}_{0,1}=
 \left(\frac{Y_{1,0}+b_{1,0}'}{X_0+ p_0'}
 ,\frac{Y_{2,0}+b_{1,0}'}{X_0+ p_0'}
 \right),\,\,\,\,\,\,
 \pmb{\rm x}_{0,2}=
  \left(\frac{Y_{1,0}+b_{1,0}''}{X_0+ p_0''}
 ,\frac{Y_{2,0}+b_{1,0}''}{X_0+ p_0''}
 \right)
 $$
 and sup-norm diameter   
 $$
  2\max_{i=0,1,2} |
   \pmb{\rm x}_0-  \pmb{\rm x}_{0,j}|_\infty\le
   2\max_{j=1,2} \,
   \max \left(
   \left|
   \frac{Y_{j,0}}{X_0}-
   \frac{Y_{j,0}+b_{j,0}'}{X_0+ p_0'}
   \right|
   ,
    \left|
   \frac{Y_{j,0}}{X_0}-
   \frac{Y_{j,0}+b_{j,0}''}{X_0+ p_0''}
   \right|,
      \left|
   \frac{Y_{j,0}}{X_0}-
   \frac{Y_{j,0}+b_{j,0}}{X_0+ p_0}
   \right|
   \right)\le
   $$
     \begin{equation}\label{suu}
     \le
     \frac{8\Delta}{X_0p_0} =
        \frac{8\Delta^3}{\gamma_1 x_0p_0^2}
        \le
           \frac{16\Delta}{\gamma_1^2p_0^2}.
  \end{equation}
  The last inequalities in (\ref{suu}) should be explained. 
  Indeed,    $$
   \left|
   \frac{Y_{j,0}}{X_0}-
   \frac{Y_{j,0}+b_{j,0}'}{X_0+ p_0'}
   \right|
   =
   \frac{|Y_{j,0}(X_0+ p_0')-  X_0(Y_{j,0}+b_{j,0}') |}{X_0 (X_0+ p_0')}        =
   \frac{|Y_{j,0}p_0'-  X_0b_{j,0}' |}{X_0 (X_0+ p_0')}    <
   $$
   $$
   <
      \frac{|Y_{j,0}p_0'-  X_0b_{j,0}' |}{X_0 ^2}
   =
   \frac{1}{X_0}
    \left|
   \frac{ p_0' Y_{j,0}}{X_0}-
   b_{j,0}'
   \right|,\,\,\,\,\,\,
   j=1,2.
   $$
   But
   \begin{equation}\label{omo}
    \left|
   \frac{ p_0' Y_{j,0}}{X_0}-
   b_{j,0}'
   \right|
   \le 
   \rho (p_0'
   \pmb{W}_0' , p_0'  \pmb{\rm x}_{0}) \le
      \rho (
  p_0' \pmb{W}_0' ,  p_0' \pmb{W}_{0})+
   p_0' \cdot
      \rho (
   \pmb{W}_0, \pmb{\rm x}_{0}) ,\,\,\,\,\,\,
   j=1,2 
   .
   \end{equation}
   For the two summands in the right hand side here we have the bound 
   $$
      \rho (
     p_0' \pmb{W}_0' ,  p_0' \pmb{W}_{0}) =
     \rho (\pmb{w}_0' ,\langle  \pmb{w}_{0}\rangle_\mathbb{R}\cap\{x_0 =p_0'\})) \le 
     2\rho (\pmb{w}_0' ,\langle  \pmb{w}_{0}\rangle_\mathbb{R})) 
     \le
     \frac{2\Delta}{|\pmb{w}_0|}
     $$ 
     and the bound of Remark 4, respectively. So we continue (\ref{omo}) with 
   $$\left|
   \frac{ p_0' Y_{j,0}}{X_0}-
   b_{j,0}'
   \right|
   \le \frac{2\Delta}{|\pmb{w}_0|} + p_0' \cdot
   \frac{2\Delta}{\gamma_1 p_0^2} \le
   \frac{4\Delta}{p_0}.
   $$
   Quite similar bounds are valid for 
   $
   \left|
   \frac{ p_0'' Y_{j,0}}{X_0}-
   b_{j,0}''
   \right|
   $
   and
   $
   \left|
   \frac{ p_0 Y_{j,0}}{X_0}-
   b_{j,0}
   \right|
   $, $ j = 1,2$. This gives the first inequality in (\ref{suu}).
 To get the last inequality in (\ref{suu}) we use  (\ref{suomi}) and  (\ref{del}).
 So we explained how to prove (\ref{suu}).
 
 So as $\pmb{W}_1 \in {\rm conv}\, (
  \pmb{\rm x}_{0}, \pmb{\rm x}_{0,0}, \pmb{\rm x}_{0,1}, \pmb{\rm x}_{0,2})$
  from (\ref{suu}) we deduce the inequality
  \begin{equation}\label{oooo}
  |\pmb{W}_1
-\pmb{\rm x}_0|
\le
\frac{16\Delta}{\gamma_1^2p_0^2}.
 \end{equation}
 This gives 
 \begin{equation}\label{gives}
 |\pmb{\rm x} -\pmb{\rm x}_0| \le
  |\pmb{W}_1
-\pmb{\rm x}_0|+
|   \pmb{W}_1
-\pmb{\rm x}|\le
{
\frac{16\Delta}{\gamma_1^2p_0^2}
+
\frac{\sqrt{2}\Delta}{\gamma_2 p_{  0}p_1} }\le 
\frac{\Delta}{2\gamma_1 p_0^2}, 
 \end{equation}
 (we used the triangle inequality, conditions (\ref{pat})  with bound $p_1 \ge p_0$ and (\ref{gammaI}))
 and we have (\ref{o}).
 
 So $\pmb{w}_0$ is a best approximation vector for  $\pmb{\rm x}$. 
 In (\ref{w2}) we have an upper bound for $\Delta_1 $ which does not depend on $\gamma$. This means that  for any 
 $\pmb{w} \in \mathbb{Z}^3\setminus 
    \langle \pmb{w}_0,\pmb{w}_1\rangle_\mathbb{R}$ we have
 $$ \rho (\pmb{w} , \langle \pmb{w}_0,\pmb{w}_1\rangle_\mathbb{R}) \ge \frac{1}{\Delta_1} \ge 
\frac{\Delta}{12p_0}.
$$
For large $\gamma_1$ the point $ \pmb{w}_0$ is essentially closer to the line $\langle \pmb{\rm x}\rangle_\mathbb{R}$
than the points $\pmb{w} = (p,b_1,b_2)\in  \mathbb{Z}^3\setminus 
    \langle \pmb{w}_0,\pmb{w}_1\rangle_\mathbb{R}$ with $ p \le p_1$.
Indeed, put $ \pmb{W} = \left( \frac{b_1}{p},\frac{b_2}{p}\right)$, then by the previous inequality and (\ref{pat},\ref{gammaI}) we see that 
$$
\sqrt{2}
| p
\pmb{\rm x} -  p\pmb{W}|_\infty\ge
| p
\pmb{\rm x} -  p\pmb{W}|
\ge \frac{1}{\Delta_1}-
 p
|\pmb{\rm x}-\pmb{W}_1|\ge
\frac{
\Delta}{12p_0} -\frac{\Delta}{\gamma_2 p_0}\ge
\frac{
\Delta}{13p_0} .
$$
At the same time 
$$
| p_0
\pmb{\rm x} -  p_0\pmb{W}_0|_\infty\le \frac{2\Delta}{\gamma_1 p_0}
$$
by (\ref{gives}) and Remark 4.
As  $\gamma_1 \ge 50$ we see that there $\pmb{\rm x}$ has no best approximations $\pmb{w}\in   \mathbb{Z}^3\setminus 
    \langle \pmb{w}_0,\pmb{w}_1\rangle_\mathbb{R}$ 
    with $ p_0 \le p \le p_1$
So  for all $\pmb{\rm x}$  satisfying  (\ref{pat}) all the best approximations between $p_0$ and $p_1$ lie in 
the two-dimensional subspace  $\langle \pmb{w}_0,\pmb{w}_1\rangle_\mathbb{R}$.
We see from  (\ref{pat})  that 
$$
|p_1\pmb{\rm x} -  p_1\pmb{W}_1|_\infty \le  \frac{\Delta}{\gamma_2 p_0}.
$$
But from the construction (\ref{norma}) and (\ref{gives}) we have
$$
\sqrt{2}
|p_0\pmb{\rm x} -  p_0\pmb{W}_0|_\infty\ge
|p_0\pmb{\rm x} -  p_0\pmb{W}_0|\ge  \frac{\Delta}{\gamma_1p_0} -  \frac{\Delta}{2\gamma_1p_0} =  \frac{\Delta}{2\gamma_1p_0}   .
$$
So
$$
|p_0\pmb{\rm x} -  p_0\pmb{W}_0|_\infty 
>
|p_1\pmb{\rm x} -  p_1\pmb{W}_1|_\infty .
$$
We have the following situation.
For any $ \pmb{\rm x}$ satisfying (\ref{pat})  vectors $\pmb{w}_0,\pmb{w}_1$ are  best approximation vectors,
and we do not have best approximation vectors $ \pmb{w} = (p,b_1,b_2) \in \mathbb{Z}^3\setminus \langle \pmb{w}_0,\pmb{w}_1\rangle_\mathbb{R}$
with $ p_0 \le p \le p_1$.
The parallelogram with vertices
$\pmb{0},
  \pmb{w}_0, \,\,\,\pmb{w}_1-\pmb{w}_0,\,\,\,
  \pmb{w}_1
 $ is a fundamental parallelogram for the lattice  $ \langle \pmb{w}_0,\pmb{w}_1\rangle_\mathbb{Z}$. So the distances from its vertices 
$  \pmb{w}_0, $ and $\pmb{w}_1-\pmb{w}_0$ to  the diagonal $ \langle \pmb{w}_1\rangle_{\mathbb{R}}$ are equal.
This means that for a point   $ \pmb{\rm x}$ which is close to $\pmb{W}_1$,  the only one possible opportunity for a vector 
$ \pmb{w} = (p,b_1,b_2) \in \mathbb{Z}^3$
with  $ p_0 \  <p < p_1$ to be a best approximation  to  $ \pmb{\rm x}$ is $ \pmb{w} =  \pmb{w}_1-\pmb{w}_0$.
Of course we cannot say that the vector $   \pmb{w}_1-\pmb{w}_0$ is a best approximation for sure. It depends on which of the  vectors
$\pmb{w}_0$ and $  \pmb{w}_1-\pmb{w}_0$ is closer to  the line spanned by the point $(1,x_1,x_2)$.

We see that   for all $\pmb{\rm x}$  satisfying  (\ref{pat}) all the best approximations with denominators between $p_0$ and $p_1$
  should be among the vectors  
$
  \pmb{w}_0, \,\,\,\pmb{w}_1-\pmb{w}_0,\,\,\,
  \pmb{w}_1
  ,$
  and everything is proved.$\Box$

     \vskip+0.3cm
  Here we should note that from (\ref{oooo}) and Remark 4  by the triangle inequality immediately follows

  \vskip+0.3cm
 {\bf Remark 7.}\, {\it
 For the rational points $\pmb{W}_0$ and $ \pmb{W}_1$ from Lemma 5 one has
 $$
 | \pmb{W}_0- \pmb{W}_1|\le \frac{3\Delta}{\gamma_1 p_0^2}.
 $$
 }
   
     \vskip+0.3cm
  {\bf 11. Proof of Theorem 2.}

  \vskip+0.3cm
  We construct a sequence of integer vectors 
  \begin{equation}\label{inn}
  \pmb{z}_\nu = (q_\nu, a_{1,\nu}, a_{2,\nu}) \in \mathbb{Z},\,\,\,\,\, \nu \in \mathbb{Z}_+
  \end{equation}
  which will be "almost" best approximation vectors to the limit point 
  \begin{equation}\label{limi}
  \pmb{\alpha} = \lim_{\nu\to \infty}
  \pmb{A}_\nu
  \end{equation}
   where
  $$
  \pmb{A}_\nu =
   \left( \frac{a_{1,\nu}}{q_\nu}, \frac{a_{2,\nu}}{q_\nu}\right) $$
   are the corresponding rational points.
  For these vectors and $ \pmb{\rm x} = (x_1,x_2)$ we consider the values
  $$
  \xi_\nu = \max_{j=1,2} | q_\nu x_j - a_{j,\nu}|,
  $$
  which of course depend on $\pmb{\rm x}$.
  
  First of all we consider  the lattice
  $$
  \Lambda_1 = 
  \langle \pmb{e}_1, \pmb{e}_2\rangle_\mathbb{Z},\,\,\,\,
  \pmb{e}_1 = (1,0,0), \,\,\,\pmb{e}_2 = (0,1,0).$$
  We put   $i_1 = 1$   
  and take
    $$
    \pmb{z}_{i_1} = \pmb{z}_1= (q_1, a_{1,1}, a_{2,1}),\,\,\,\pmb{z}_{i_1+1} = \pmb{z}_2=(q_2, a_{1,2}, a_{2,2}) 
    $$
    to be a basis of $\Lambda_1$  in such a way that $ q_2>q_1$ and all  the  conditions  (\ref{poo}, \ref{x1})
 of Lemma 4 are satisfied  
 for $ \pmb{v}_0 = \pmb{z}_1,  \pmb{v}_1 = \pmb{z}_2$.
   (In particular, the condition (\ref{x1}) is satisfied if the angle between the basis  vectors  $  \pmb{z}_1,  \pmb{z}_2$ is small.)
      We take
   $$
   \gamma =  \max ( 400, q_2/q_1),
   $$
so $ q_2 \le \gamma q_1$.
   Now we 
  define  vectors (\ref{inn})  by inductive procedure.
  Let 
 vectors (\ref{inn}) be defined up to $ \pmb{z}_{i_t+1}$ and  the following conditions are valid

   \vskip+0.3cm 
 \noindent 
{\rm ({\bf A})} 
 two last vectors 
 $ \pmb{v}_0 = \pmb{z}_{i_t}, \pmb{v}_1 =  \pmb{z}_{i_t+1}$  satisfy  all the conditions   (\ref{poo}, \ref{x1}) of Lemma 4
  where  
  $ p_0 = q_{i_t}, p_1 = q_{i_t+1}$
  and
  $ \Delta = \Delta _t$ is the
   fundamental volume of two-dimensional lattice $ \Lambda_t =   \langle \pmb{z}_{i_t}, \pmb{z}_{i_t+1} \rangle_\mathbb{Z}$,
   moreover 
   \begin{equation}\label{aq}
      |\pmb{V}_{0}-\pmb{V}_{1}|_\infty=
   |\pmb{A}_{i_t}-\pmb{A}_{i_t+1}|_\infty
   \le \frac{\Delta_t}{30\gamma^2 q_{i_t}^2}
   ;
 \end{equation}

   \vskip+0.3cm 
 \noindent 
{\rm ({\bf B})}  for all $ \pmb{\rm x} = (x_1,x_2)$  satisfying 
\begin{equation}\label{saa}
  | \pmb{\rm x} -\pmb{A}_{i_t}|_\infty
\le \frac{\Delta_t}{24\gamma^2 q_{i_t}^2}
\end{equation}
all the best approximation vectors $ \pmb{z} = (q, a_1,a_2)$  with $ q_1\le q \le q_{i_t }$ are among the vectors from
the sequence 
\begin{equation}\label{ssq}
 \pmb{z}_1,\pmb{z}_2,...,\pmb{z}_{i_{t} };
 \end{equation}

   \vskip+0.3cm

   \noindent 
{\rm ({\bf C})}
among every two consecutive vectors 
$\pmb{z}_\nu,\pmb{z}_{\nu+1}$ from  (\ref{ssq}) at least one vector is a best approximation vector for every 
$ \pmb{\rm x}  $  satisfying (\ref{saa});

   \vskip+0.3cm 
 \noindent 
{\rm ({\bf D})}  for all $\pmb{\rm x}$ satisfying (\ref{saa}) and for every $ \nu \le i_t-1$ one has 
 $
 \frac{\xi_\nu}{\xi_{\nu-1}} \ge \frac{1}{16\sqrt{6}(50\gamma^2+2)}
 $.

   \vskip+0.3cm

   When vectors (\ref{inn}) will be defined, the limit point  (\ref{limi}) will 
   satisfy
   \begin{equation}\label{saar}
  | \pmb{\alpha} -\pmb{A}_{i_t}|_\infty
\le \frac{\Delta_t}{24\gamma^2 q_{i_t}^2} \,\,\,\forall t \in \mathbb{Z}_+,
\end{equation}
as the inequality (\ref{saa}) for $t+1$ leads to  the inequality (\ref{saa}) for $t$.
This limit vector $\pmb{\alpha}$ will
   be just the vector we need for Theorem 2.

   \vskip+0.3cm 
   Here we should note that  for $ t=1$ the conditions  {\rm ({\bf B})}   is satisfied automatically as
  $ \pmb{z}_1$ is a best approximation vector  for all  $ \pmb{\rm x}$ satisfying (\ref{saa}). At the same time
   for $t=1$ conditions  {\rm ({\bf C})}  and {\rm ({\bf D})} are empty, because we have only one vector $\pmb{z}_1$.
   \vskip+0.3cm 
 
 Now we explain how to construct next vectors 
 \begin{equation}\label{vectors}
 \pmb{z}_\nu,\,\,\,\,\,\,\, i_t+2\le \nu \le i_{t+1}+1.
 \end{equation}
  satisfying conditions 
   {\rm ({\bf A})}, {\rm ({\bf B})}, {\rm ({\bf C})}, {\rm ({\bf D})} of the next step.
We start with the explanation of  the construction  and then we will verify the conditions 
   {\rm ({\bf A})}, {\rm ({\bf B})}, {\rm ({\bf C})}, {\rm ({\bf D})} .
   
    First of all we   apply Lemma 4
   with 
   $$
   \pmb{v}_0 = \pmb{z}_{i_t},\,\,\,\,
   \pmb{v}_1= \pmb{z}_{i_t+1}
   $$
   and
   take vectors 
   $$
   \pmb{z}_{i_t+\nu }  =\pmb{v}_\nu, \,\,\,\, 2\le \nu \le k_t
   $$
   where $\pmb{v}_\nu$ are defined in  (\ref{x2}).
   We take $k= k_t$ large enough to satisfy (\ref{u}) as well as the inequalities
   \begin{equation}\label{qqu}
   q_{i_t+k_t} \ge \gamma \Delta_{t}^2
   \end{equation}
   and
     \begin{equation}\label{qqu1}
   q_{i_t+k_t} \ge \gamma    q_{i_t},
   \end{equation}
   
     We define
   $$
   i_{t+1} = i_t+k_t +2,
   $$
   so
   $$
   i_{t+1}-2 =  i_t+k_t.
   $$
   Then for 
   $$
   \gamma_1 = \gamma,\,\,\,\gamma_2 = \gamma^2
   $$
   and
   vectors
   $$
   \pmb{w}_0' = \pmb{z}_{i_{t+1}-4},\,\,\,\,
 \pmb{w}_0''= \pmb{z}_{i_{t+1}-3}
 $$ we apply Lemma 5.
 Of course we have 
 $$ \pmb{w}_0 = \pmb{w}_0'+\pmb{w}_0''=  \pmb{z}_{i_{t+1}-2}.
 $$
 We need to check the condition
 (\ref{del}) and the condition on $\pmb{\rm x}$ satisfying (\ref{o}).
 But (\ref{del}) follows from (\ref{qqu}). As for 
  the condition on $\pmb{\rm x}$ we will check it right now.
  In our situation $ p_0 = q_{i_{t+1}-2}$ and (\ref{o}) means that
    $$
    |\pmb{\rm  x} -   \pmb{\rm  x}_0|_\infty \le \frac{\Delta_t}{\gamma  q_{i_{t+1}-2}^2}
  .
  $$
  Remark 4  with $ \pmb{W}_0 = \pmb{A}_{i_{t+1}-2}$ gives
  $$
  |\pmb{\rm  x}_0 -   \pmb{W}_0|_\infty=
    |\pmb{\rm  x}_0 -   \pmb{A}_{i_{t+1}-2}|_\infty
  \le
  \frac{2\Delta_t}{\gamma q_{i_{t+1}-2}^2}.
  $$
  So by the triangle inequality
  $$ 
    |\pmb{\rm  x} -   \pmb{A}_{i_{t+1}-2}|_\infty
  \le
  \frac{3\Delta_t}{\gamma q_{i_{t+1}-2}^2}\le
  \frac{\Delta_t}{100 q_{i_{t+1}-2}^2}.
  $$
  So  $\pmb{\rm x}$ satisfies (\ref{x3})  and the condition on $\pmb{\rm x}$ follows from the conclusion of Lemma 4,
  as in  both sequences (\ref{x5}) and (\ref{x51}) the last vector is
  $\pmb{v}_k = \pmb{w}_0 = \pmb{z}_{i_{k+1}-2}$.
  We verified the possibility of application of Lemma 5.
  Lemma 5  gives us the vector
 $$
 \pmb{z}_{i_{t+1}} =  \pmb{w}_1.
 $$
 Then we define
 $$
  \pmb{z}_{i_{t+1}-1} =  \pmb{w}_1 - \pmb{w}_0.
 $$
 Now we should define  $
 \pmb{z}_{i_{t+1}+1} .
 $ First of all we define the next two-dimensional lattice
 $\Lambda_{t+1} = 
 \langle  \pmb{z}_{i_{t+1}-1},  \pmb{z}_{i_{t+1}}\rangle_\mathbb{Z}$ with fundamental volume
 $\Delta_{t+1}$.
 Then we define
 \begin{equation}\label{30}
  \pmb{z}_{i_{t+1}+1} =  \pmb{z}_{i_{t+1}-1}+a
   \pmb{z}_{i_{t+1}}, \,\,\,\,\,
   \text{where}\,\,\,\,\,\,\,
   a= [50\gamma^2 ] +1.
   \end{equation}
   It is clear that 
    $\Lambda_{t+1} = 
 \langle  \pmb{z}_{i_{t+1}},  \pmb{z}_{i_{t+1}+1}\rangle_\mathbb{Z}$.
 
 So all the vectors  (\ref{vectors})
 are defined and we must check the conditions 
   {\rm ({\bf A})}, {\rm ({\bf B})}, {\rm ({\bf C})}, {\rm ({\bf D})} of the new inductive step.

 \vskip+0.3cm
 Condition 
   {\rm ({\bf A})}
 is satisfied because  of
 $$
 |
 a_{j,i_{t+1}+1}  q_ {i_{t+1}}-
 a_{j, i_{t+1}}q_{i_{t+1}+1}|\le
 \Delta_{t+1},
 $$
 and for (\ref{aq}) with $t$ replaced by $t+1$ we have
 \begin{equation}\label{l0}
     |\pmb{A}_{i_{t+1}} - \pmb{A}_{i_{t+1}+1}|_\infty = 
         \max_{j=1,2} \left| \frac{a_{j,i_{t+1}+1}}{q_{i_{t+1}+1}}- \frac{a_{j, i_{t+1}}}{q_ {i_{t+1}}}\right| 
\le \frac{\Delta_{t+1}}{  q_{i_{t+1}}q_{i_{t+1}+1}}\le
\frac{\Delta_{t+1}}{ A q_{i_{t+1}}^2} \le 
    \frac{\Delta_{t+1}}{ 50\gamma^2  q_{i_{t+1}}^2} 
         .
    \end{equation}
    Let us check the conditions of Lemma 4.
 Inequality  (\ref{poo}) is clear. As for (\ref{x1}), we should show that 
 \begin{equation}\label{l1}
    |\pmb{A}_{i_{t+1}} - \pmb{A}_{i_{t+1}+1}|_\infty\le \frac{1}{2q_{i_{t+1}}\Delta_{t+1}}
 \end{equation}
 and
  \begin{equation}\label{l29}
   |\pmb{A}_{i_{t+1}} - \pmb{A}_{i_{t+1}+1}|_\infty\le \frac{\delta ( \pmb{z}_{i_{t+1}} )}{2}.
 \end{equation}
 To get  (\ref{l1}) we use (\ref{dee}) for the best approximation vector $\pmb{z}_{i_{t+1}} $
with $ \Delta_{t+1}$ instead of $\Delta_2$. Then 
$$
q_{i_{t+1}}\ge (K\Delta_{t+1})^2\ge 
 \frac{\Delta_{t+1}^2}{\gamma^2},
 $$
 and this deduces (\ref{l1}) from ({\ref{l0}).

From  condition ({\bf v}) of Lemma 5 we see that  
$ \pmb{w}_1 = \pmb{z}_{i_{t+1}} $  will be a best approximation vector for all $\pmb{\rm x}$ satisfying the condition  (\ref{pat}). So we have
$$ 
\frac{\Delta_t}{\gamma^2 q_{i_{t+1}-2}q_{i_{t+1}}} \le \delta (\pmb{z}_{i_{t+1}} ),
$$
as $\pmb{w}_0 =\pmb{z}_{i_{t+1}} $ is always a best approximation vector under the assumption (\ref{pat}).
The last inequality together with (\ref{l0}) and Remark 5 
(where $ \Delta = \Delta_t, \Delta_1 = \Delta_{t+1}, p_0 = q_{i_{t+1}-2}, p_1 = q_{i_{t+1}}$)
gives
$$
     |\pmb{A}_{i_{t+1}} - \pmb{A}_{i_{t+1}+1}|_\infty\le\frac{\Delta_{t+1}}{50\gamma^2 q_{i_{t+1}}^2}
     \le
     \frac{\Delta_{t}}{2\gamma^2 q_{i_{t+1}-2} q_{i_{t+1}}} \le \frac{\delta (\pmb{z}_{i_{t+1}} )}{2},
     $$
     and this is just (\ref{l29}). 
     
     So condition  {\rm ({\bf A})} is satisfied.
 
  \vskip+0.3cm
  
  Now we verify conditions     {\rm ({\bf B})}  and   {\rm ({\bf C})}.
  Suppose that $ \pmb{\rm x}$ satisfies  (\ref{saa}) for the next step, that is 
  \begin{equation}\label{saa1}
  | \pmb{\rm x} -\pmb{A}_{i_{t+1}}|_\infty
\le \frac{\Delta_{t+1}
}{24\gamma^2 q_{i_{t+1}}^2}
\end{equation}
From (\ref{saa1})  and Remark 5 we see that
$$
  | \pmb{\rm x} -\pmb{A}_{i_{t+1}}|_\infty
\le \frac{\Delta_{t}
}{\gamma^2  q_{i_{t+1}-2}q_{i_{t+1}}}.
  $$
  So  by Lemma 5 either 
  $$
  \pmb{z}_{i_{t+1}-2} ,\pmb{z}_{i_{t+1}-1} ,\pmb{z}_{i_{t+1}} 
  $$
  or
  $$
  \pmb{z}_{i_{t+1}-2}  ,\pmb{z}_{i_{t+1}} 
  $$
  are successive best approximations to $\pmb{\rm x}$.

  Then
  from Remark 7
  (with $\pmb{W}_0 = \pmb{A}_{i_{t+1}-2},
  \pmb{W}_1= \pmb{A}_{i_{t+1}}, \Delta= \Delta_t, \gamma_1= \gamma, p_0 = q_{i_{t+1}-2}^2$)
  we have
  $$
    | \pmb{A}_{i_{t+1}-2} -\pmb{A}_{i_{t+1}}|_\infty
    \le
       \frac{3\Delta_t}{400 q_{i_{t+1}-2}^2}.
       $$
       This inequality together with 
  (\ref{saa1}) leads to
  $$
    | \pmb{\rm x} -\pmb{A}_{i_{t+1}-2}|_\infty 
    \le
      | \pmb{\rm x} -\pmb{A}_{i_{t+1}}|_\infty
      +
        | \pmb{A}_{i_{t+1}-2} -\pmb{A}_{i_{t+1}}|_\infty
        \le
        \frac{\Delta_t}{100 q_{i_{t+1}-2}^2}.
        $$
      So by Lemma 4 we see that either
      $$
       \pmb{z}_{i_{t}} ,\pmb{z}_{i_{t}+1} ,\pmb{z}_{i_{t+1}-4} ,
        \pmb{z}_{i_{t+1}-3} ,\pmb{z}_{i_{t+1}-2}  
      $$
      or
       $$
       \pmb{z}_{i_{t}} ,\pmb{z}_{i_{t}+1} ,\pmb{z}_{i_{t+1}-4} 
       ,\pmb{z}_{i_{t+1}-2}  
      $$
      is the sequence of successive best approximations to  $\pmb{\rm x}$.

      Again from (\ref{saa1}) 
      and Remark 6 
      ($ q_{i_{t}}\le p_0 =q_{i_{t+1}-2}, p_1 = q_{i_{t+1}}$) which now states that
      $ q_{i_{t+1}}\ge \frac{\gamma}{2} q_{i_{t+1}-2}$
      we deduce
      $$
        | \pmb{\rm x} -\pmb{A}_{i_{t+1}}|_\infty\le
        \frac{\Delta_t}{24 \gamma^2 q_{i_{t+1}}^2}\le
           \frac{\Delta_t}{96 \gamma^4 q_{i_{t}}^2} 
      .$$
      Then, by Remark 7 ($
      \pmb{W}_0 =  \pmb{A}_{i_{t+1}-2}, \pmb{W}_1 =  \pmb{A}_{i_{t+1}}
      p_0 =  q_{i_{t+1}-2} = q_{i_t + k_t}
      $)
       and (\ref{qqu1}) we see that 
      $$
       | \pmb{A}_{i_{t+1}}-\pmb{A}_{i_{t+1}-2}|_\infty\le
           \frac{3\Delta_t}{\gamma_1 p_0^2} =
             \frac{3\Delta_t}{\gamma q_{i_{t+1}-2}^2}
       \le
           \frac{3\Delta_t}{ \gamma^3q_{i_{t}}^2} 
             \le
           \frac{\Delta_t}{ 100\gamma^2q_{i_{t}}^2} 
           .
      $$
     In the notation of Lemma 4
      we have $\pmb{V}_0 = \pmb{A}_{i_{t}}, \pmb{V}_1= \pmb{A}_{i_{t}+1},
      \pmb{V}_k = \pmb{A}_{i_{t+1}-2}$.
   So
      $$
           | \pmb{A}_{i_{t+1}-2}-\pmb{A}_{i_{t}}|_\infty
           =|\pmb{V}_k-\pmb{V}_0|_\infty 
           \le
           |\pmb{V}_1-\pmb{V}_0|_\infty =
            | \pmb{A}_{i_{t}+1}-\pmb{A}_{i_{t}}|_\infty
              \le \frac{\Delta_t}{30\gamma^2 q_{i_t}^2}
            ,
      $$
      by (\ref{aq}) from condition  {\rm ({\bf A})}.
      So last three inequalities lead to 
      $$
       | \pmb{\rm x} -\pmb{A}_{i_{t}}|_\infty
\le 
  | \pmb{\rm x} -\pmb{A}_{i_{t+1}}|_\infty+
    | \pmb{A}_{i_{t+1}}-\pmb{A}_{i_{t+1}-2}|_\infty
    +
       | \pmb{A}_{i_{t+1}-2}-\pmb{A}_{i_{t}}|_\infty
\le
\frac{\Delta_{t}
}{24\gamma^2 q_{i_{t}}^2},
$$
and by inductive assumption we have the required properties for all the best approximations
$ \pmb{z} = (q, a_1,a_2)$  with $ q_1\le q \le q_{i_t }$.

By the way,  we see that condition (\ref{saa}) for  $(t+1)$-th step ensures condition  (\ref{saa}) for $t$-th step,
and we proved the inequality (\ref{saar}).

We see that we have established conditions  {\rm ({\bf B})} and  {\rm ({\bf C})} for all the best appproximations 
$ \pmb{z} = (q, a_1,a_2)$  in the range $ q_1\le q \le q_{i_{t+1 }}$.
      
      \vskip+0.3cm
      
      Let us verify condition {\rm ({\bf D})} 
      for $ i_t \le \nu \le i_{t+1}-1$.
      We consider the cases
      
      1)  $\nu = i_t$,
      
      2) $ i_t < \nu \le i_{t+1}-3$,
      
      3) $\nu =  i_{t+1}-2$,
      
       4) $\nu =  i_{t+1}-1$
       
       \noindent
       separately.
       
         \vskip+0.3cm
         
         1) First of all we need lower bound for the approximation
      $
        \xi_{i_t} = q_{i_t}  | \pmb{\alpha} -\pmb{A}_{i_t}|_\infty
 .
$
We use the notation of Lemma 4 with 
$$
\frak{y}_0 =  q_{i_t}  \underline{\pmb{x}} -\underline{\pmb{z}}_{i_t},\,\,\,\,\,\,
  \xi_{i_t} = |\underline{\frak{y}}_0 |_\infty.
  $$
  By (\ref{30}) of the previous inductive step we have
    $$
   q_{i_t+1}\le (50\gamma^2+2) q_{i_t}
   .$$
   Remark 5 with $ p_0 =q_{i_t}, p_1 = q_{i_t+1}, \Delta = \Delta_{t-1}, \Delta_1 = \Delta_t$ for Lemma 5 applied on the previous inductive step gives
      $$
   \Delta_t \ge \frac{q_{i_t}}{8q_{i_t-2}}\Delta_{t-1}
   $$
  Now
  from (\ref{l2}) with $ i=0, \Delta= \Delta_{t}, $
 we get
 
\begin{equation}\label{xi1}
\xi_{i_t} \ge \frac{\Delta_t}{2\sqrt{6}q_{i_t+1}}\ge
\frac{\Delta_{t-1}}{16\sqrt{6}(50\gamma^2 +2)q_{i_t-2}}
.
\end{equation}
Form (\ref{0000})  with $ \nu =i_t-2, \Delta_2 = \Delta_{t-1}$ we see that
\begin{equation}\label{xi2}
\xi_{i_{t-2}} \le \frac{\Delta_{t-1}}{ q_{i_t-1}} \le \frac{\Delta_{t-1}}{ q_{i_t-2}}.
\end{equation}
Points $\pmb{0}, \pmb{z}_{i_t-2}, \pmb{z}_{i_{t}-1},\pmb{z}_{i_t}$ form a parallelogram and so
$$
\xi_{i_t-1}= \xi_{i_t-2}-\xi_{i_t} < \xi_{i_t-2}.
$$
 Now (\ref{xi1},\ref{xi2}) give us
 $$
      \frac{\xi_{i_t}}{\xi_{i_{t-1}-1}} \ge \frac{\xi_{i_t}}{\xi_{i_{t-1}-2}}\ge\frac{1}{16\sqrt{6}(50\gamma^2 +2)},
      $$
      and this is what we need.
          \vskip+0.3cm 
      
     2) For $\nu$ from the interval $ i_t < \nu \le i_{t+1}-3$ from (\ref{laci}) of Lemma 4  follows
      $$
      \frac{\xi_{\nu}}{\xi_{\nu-1}} \ge \frac{1}{4}.
      $$

       \vskip+0.3cm
       
     3) Let $\pmb{\rm x}_0 $ be the point form Lemma 5 applied on $(t+1)$-th step.
     In the notation of Lemma 5 we have $ \pmb{W}_0 = \pmb{A}_{i_{t+1}-2}, \pmb{W}_1 = \pmb{A}_{i_{t+1}},
     p_0 = q_{i_{t+1}-2}, p_1 = q_{i_{t+1}}, \Delta = \Delta_t, \Delta_1 = \Delta_{t+1}$. Then
     \begin{equation}\label{mm1}
     \xi_{i_{t+1}-2} \ge
     q_{i_{t+1}-2}|\pmb{A}_{i_t{+1}-2} -\pmb{\rm x}_0 |_\infty -
     q_{i_{t+1}-2}|\pmb{x} -\pmb{\rm x}_0 |_\infty.
     \end{equation}
     But from the construction  (\ref{norma}) we have
         \begin{equation}\label{mm2}
     q_{i_{t+1}-2}|\pmb{A}_{i_t{+1}-2} -\pmb{\rm x}_0 |_\infty 
     =
     p_0 |\pmb{W}_0 - \pmb{\rm x}_0 |_\infty  \ge
     \frac{\Delta_t}{2\gamma q_{i_{t+1}-2}}.
     \end{equation}
     Then, 
        \begin{equation}\label{mm3}
    q_{i_{t+1}-2}|\pmb{x} -\pmb{\rm x}_0 |_\infty
    \le
     q_{i_{t+1}-2}| \pmb{A}_{i_{t+1}} -\pmb{\rm x}_0 |_\infty
     +
     p_0 |\pmb{\rm x} -  \pmb{A}_{i_{t+1}} |_\infty  \le
     \frac{17\Delta_t}{\gamma^2q_{i_{t+1}-2}}
     \end{equation}
     (we use inequalities (\ref{oooo}) for the first summand and    $(t+1)$-th step of (\ref{saa}), Remark 5 for the second summand).
     Now (\ref{mm1},\ref{mm2},\ref{mm3}) gives
     $$
     \xi_{i_{t+1}-2} \ge \frac{\Delta_t}{4\gamma q_{i_{t+1}-2}}.
     $$
     Together with (\ref{0000}) for
     $ \nu =  {i_{t+1}-3}$
     this gives
     $$
     \frac{ \xi_{i_{t+1}-2} }{ \xi_{i_{t+1}-3} }\ge\frac{1}{4\gamma}.
     $$
     
           \vskip+0.3cm
       
     4) As in the case 1) the points
     $\pmb{0}, \pmb{z}_{i_{t+1}-2}, \pmb{z}_{i_{t+1}-1},\pmb{z}_{i_{t+1}}$ form a parallelogram and so
$$
\xi_{i_{t+1}-1}= \xi_{i_{t+1}-2}-\xi_{i_{t+1}}.
$$
As  $\xi_{i_{t+1}}$ 
     is much smaller than
     $\xi_{i_{t+1}-2}$ 
     we immediately have
    $$
     \frac{ \xi_{i_{t+1}-1} }{ \xi_{i_{t+1}-2} }=
     1-
     \frac{ \xi_{i_{t+1}} }{ \xi_{i_{t+1}-2} }
     \ge\frac{1}{2}.
     $$
            \vskip+0.3cm
      We see that condition {\rm ({\bf D})}  is valid in the range  $ i_t < \nu \le i_{t+1}$.
            \vskip+0.3cm
            
      Now we have constructed the vectors (\ref{inn}) satisfying the conditions 
         {\rm ({\bf A})}, {\rm ({\bf B})}, {\rm ({\bf C})}, {\rm ({\bf D})}   for every $t$ and Theorem 2 follows.$\Box$.
      
      \vskip+0.3cm
      
     {\bf Acknowledgements.}
     
       \vskip+0.3cm
     
     The authors  thank the anonymous referee for careful reading of the manuscript and for important suggestions.
     
       The second named  is a winner of the   “Leader” contest conducted by Theoretical Physics and Mathematics Advancement Foundation “BASIS” and would like to thank the foundation and jury.

\vskip+1cm

Renat Akhunzhanov,

Astrakhan State University,

e-mail:  akhunzha@mail.ru

\vskip+1cm

Nikolay Moshchevitin

Moscow Center of Fundamental and Applied Mathematics

and

Steklov Mathematical Institute

e-mail: Moshchevitin@gmail.com

\end{document}